\documentclass{gtart}

\def\ifplaintex{\expandafter\ifx\csname documentclass\endcsname\relax}

\expandafter\ifx\csname epsfbox\endcsname\relax\input epsf\fi

\ifplaintex 
\else
\fi

\def\gt{{\mathsurround=0pt\it $\cal G\mskip-2mu$eometry \&\ 
$\cal T\!\!$opology}}        %  journal title in recommended style

\def\gtp{{\mathsurround=0pt\it $\cal G\mskip-2mu$eometry \&\ 
$\cal T\!\!$opology $\cal P\!$ublications}}  % GT publications

\def\lognumber#1{\def\thelognumber{#1}}
\def\volumenumber#1{\def\thevolumenumber{#1}}
\def\papernumber#1{\def\thepapernumber{#1}}
\def\volumeyear#1{\def\thevolumeyear{#1}}

\def\pagenumbers#1#2{\def\startpage{#1}\def\finishpage{#2}}
\def\published#1{\def\publishdate{#1}}
\def\proposed#1{\def\theproposer{#1}}
\def\seconded#1{\def\theseconders{#1}}
\def\received#1{\def\receiveddate{#1}}
\def\revised#1{\def\reviseddate{#1}}
\def\accepted#1{\def\accepteddate{#1}}

\long\def\asciiabstract#1{\long\def\theasciiabstract{#1}}

\let\\\par\let\thevolumenumber\relax\let\thepapernumber\relax
\let\thevolumeyear\relax\let\thesamplenumber\relax\let\startpage\relax
\let\finishpage\relax\let\publishdate\relax\let\receiveddate\relax
\let\reviseddate\relax\let\accepteddate\relax\let\theasciititle\relax
\let\theasciiauthors\relax
\let\theasciiabstract\relax
\let\theasciiemail\relax\let\theshortauthors\relax\let\theshorttitle\relax

\long\def\maketitlep{   % start of definition of \maketitlep

\count0=\startpage

\gt\hfill      %   Journal title (top left) 
\hbox to 77pt{\vbox to 0pt{\vglue -15pt\epsfbox{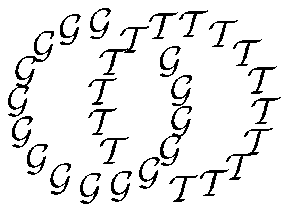}\vss}\hss}
\break
{\small\ifx\thesamplenumber\relax % sample?  
Volume \else Sample
\fi\thevolumenumber\ (\thevolumeyear)
\startpage--\finishpage\nl
Published: \publishdate}
\vglue 0.5truein plus 0.4fil minus 0.1truein

{\parskip=0pt\leftskip 0pt plus 1fil\def\\{\par\smallskip}{\ifplaintex\large
\else\Large\fi\bf\thetitle}\par\medskip}   

\vglue 0pt plus 0.1fil 

{\parskip=0pt\leftskip 0pt plus 1fil\def\\{\par}{\sc\theauthors}
\par\medskip}

\vglue 0pt plus 0.1fil 

{\small\parskip=0pt\let\newline\\
{\leftskip 0pt plus 1fil\def\\{\par}{\sl\theaddress}\par}
\expandafter\ifx\theemail\relax    % email address?
\relax\else\vglue 5pt plus 0.02fil minus 2pt\def\\{\stdspace{\rm 
and}\stdspace} 
\cl{Email:\stdspace\tt\theemail}\fi
\ifx\theurl\relax                  % URL given?
\relax\else\vglue 5pt plus 0.02fil minus 2pt\def\\{\stdspace{\rm 
and}\stdspace}
\cl{URL:\stdspace\tt\theurl}\fi\par}

\vglue 7pt plus 0.3fil minus 3pt

{\bf Abstract}
\vglue 5pt plus 0.1fil minus 2pt

\theabstract

\vglue 7pt plus 0.3fil minus 3pt

{\bf AMS Classification numbers}\quad Primary:\quad \theprimaryclass

Secondary:\quad \thesecondaryclass

\vglue 5pt plus 0.3fil minus 2pt

{\bf Keywords:}\quad \thekeywords

\vglue 10pt plus 0.5fil minus 5pt

{\small  Proposed: \theproposer\hfill Received: \receiveddate\nl
Seconded: \theseconders\hfill 
\ifx\reviseddate\relax                         % paper revised?
Accepted: \accepteddate                        % no
\else
Revised: \reviseddate                          % yes
\fi}
\eject
}       %  end of definition of \maketitlep

\let\maketitlepage\maketitlep
\let\maketitle\maketitlepage

\font\phead=cmsl9 scaled 950
\font=cmsl9 scaled 1050
\font\pnum=cmbx10 scaled 913
\font\lnum=cmbx10 
\font\pfoot=cmsl9 scaled 950
\font=cmsl9 scaled 1050
\ifplaintex
\headline{\vbox to 0pt{\vskip -4.5mm\line{\small\phead\ifnum
\count0=\startpage ISSN 1364-0380 (on line)
1465-3060 (printed) \hfill {\pnum\folio}\else\ifodd\count0\def\\{ }% 
\ifx\theshorttitle\relax\thetitle\else\theshorttitle\fi\hfill{\pnum\folio}
\else\def\\{ and }{\pnum\folio}\hfill\ifx\theshortauthors\relax\theauthors
\else\theshortauthors\fi\fi\fi}\vss}}
\footline{\vbox to 0pt{\vglue 0mm\line{\small\pfoot\ifnum\count0=\startpage
\copyright\ \gtp\hfill\else
\gt, Volume \thevolumenumber\ (\thevolumeyear)\hfill\fi}\vss
}}
\else
\makeatletter
\def\@oddhead{{\small\lhead\ifnum\count0=\startpage ISSN 1364-0380 (on line)
1465-3060 (printed) \hfill {\lnum\number\count0}\else\ifodd\count0
\def\\{ }\ifx\theshorttitle\relax \thetitle \else\theshorttitle\fi\hfill
{\lnum\number\count0}\else\def\\{ and }{\lnum\number\count0}
\hfill\ifx\theshortauthors\relax 
\theauthors\else\theshortauthors\fi\fi\fi}}\def\@evenhead{@oddhead}
\def\@oddfoot{\small\lfoot\ifnum\count0=\startpage\copyright\ \gtp\hfill\else
\gt, Volume \thevolumenumber\ (\thevolumeyear)\hfill\fi}
\def\@evenfoot{@oddfoot}
\makeatother
\fi

\newwrite\gtoutfile
\long\gdef\makeheadfile{  %%% start of definition of \makeheadfile
{\def\\{, }\def\s{ }
\immediate\openout\gtoutfile head.xxx
\immediate\write\gtoutfile{To: math@arxiv.org}
\immediate\write\gtoutfile{Subject: put OR rep NNNNN:pppp}
\immediate\write\gtoutfile{--text follows this line--}
\immediate\write\gtoutfile{Proxy-for: \ifx\theasciiauthors\relax
\theauthors\else\theasciiauthors\fi\s<\ifx\theasciiemail\relax\theemail\else\theasciiemail\fi>}
\immediate\write\gtoutfile{\noexpand\\}
\immediate\write\gtoutfile{Authors: \ifx\theasciiauthors\relax
\theauthors\else\theasciiauthors\fi}
{\def\\{ }\immediate\write\gtoutfile{Title: \ifx\theasciititle\relax
\thetitle\else\theasciititle\fi}}
\immediate\write\gtoutfile{Subj-class: GT or GR or SG or ...}
\immediate\write\gtoutfile{MSC-class: \theprimaryclass\ifx\thesecondaryclass\relax\else, \thesecondaryclass\fi}
\immediate\write\gtoutfile{Journal-ref: Geom. Topol. \thevolumenumber\s
(\thevolumeyear) \startpage-\finishpage}
\immediate\write\gtoutfile{Comments: Published in Geometry and Topology at}
\immediate\write\gtoutfile{    http://www.maths.warwick.ac.uk/gt/GTVol\thevolumenumber/paper\thepapernumber.abs.html}
\immediate\write\gtoutfile{\noexpand\\}
\immediate\write\gtoutfile{}
\ifx\theasciiabstract\relax
\immediate\write\gtoutfile{\theabstract}\else
\immediate\write\gtoutfile{\theasciiabstract}\fi
\immediate\write\gtoutfile{}
\immediate\write\gtoutfile{\noexpand\\}
\immediate\write\gtoutfile{}
\immediate\closeout\gtoutfile}}  %%% end of definition of \makeheadfile

\def\maketitlepage{\maketitlep\makeheadfile}
\let\maketitle\maketitlepage

\lognumber{156}
\volumenumber{5}\papernumber{26}\volumeyear{2001}
\pagenumbers{831}{883}
\received{18 December 2000}
\accepted{30 October 2001}
\proposed{Ronald Stern}
\seconded{Robion Kirby, Cameron Gordon}
\revised{17 November 2001}
\published{21 November 2001}

\usepackage{amssymb, amscd, amsmath, epsf}

\def\int{\mathop{\rm int}}

\def\si{\sigma}

\def\ep{\epsilon}

\def\zz{{\bf Z}}
\def\qq{{\bf Q}}
\def\cc{{\bf C}}
\def\rr{{\bf R}}

\def\calc{{\cal C}}
\def\calg{{\cal G}}
\def\calm{{\cal M}}
\def\cals{\cal S}

\def\mapright#1{\smash{ \mathop{\longrightarrow}\limits^{#1}}}

\newtheorem{thm}{Theorem}[section]
\newtheorem{lemma}[thm]{Lemma}

\newtheorem{cor}[thm]{Corollary}
\newtheorem{defn}[thm]{Definition}

\begin{document}
\title{Concordance and Mutation}
\author{P Kirk and C Livingston}
\address{Department of Mathematics\\Indiana 
University, Bloomington, IN 47405, USA}
\email{pkirk@indiana.edu,
livingst@indiana.edu}

\begin{abstract}
{\small
We provide a framework for studying the interplay between concordance
and positive mutation and identify some of the basic structures
relating the two.

The fundamental result in understanding knot concordance is the
structure theorem proved by Levine: for $n>1$ there is an isomorphism
$\phi$ from the concordance group $\calc_n$ of knotted $(2n-1)$--spheres
in
$S^{2n+1}$ to an algebraically defined group $\calg_{\pm}$; furthermore,
$\calg_{\pm}$ is isomorphic to the infinite direct sum $\zz^\infty
\oplus \zz_2^\infty \oplus \zz_4^\infty$.  It was a startling 
consequence of the work of Casson and Gordon that in the classical
case  the kernel of
$\phi$ on $\calc_1$ is infinitely generated. Beyond this, little has been
discovered about the pair $(\calc_1,\phi)$.  

In this paper we present a new approach to studying
$\calc_1$ by introducing a group,
$\calm$, defined as the quotient of the set of knots by the equivalence
relation generated by concordance and positive mutation, with group
operation induced by connected sum.  We prove there is a factorization
of
$\phi$, $\calc_1 \mapright{\phi_1} \calm \mapright{\phi_2} \calg_-$.  
Our main result is that both maps have infinitely generated kernels.

Among geometric constructions on classical knots, the most subtle is
positive mutation.  Positive mutants are indistinguishable using
classical abelian knot invariants as well as by such modern invariants
as the Jones, Homfly or Kauffman polynomials.  Distinguishing positive
mutants up to concordance is a far more difficult problem; only one
example has been known until now.  The results in this paper provide,
among other results, the first infinite families of knots that are
distinct from their positive mutants, even up to concordance.}
 \end{abstract}

\asciiabstract{We provide a framework for studying the interplay
between concordance and positive mutation and identify some of the
basic structures relating the two.

The fundamental result in understanding knot concordance is the
structure theorem proved by Levine: for n>1 there is an isomorphism
phi from the concordance group C_n of knotted (2n-1)-spheres in
S^{2n+1} to an algebraically defined group G_{+-};
furthermore, G__{+-} is isomorphic to the infinite direct sum
Z^infty direct sum Z_2^infty direct sum Z_4^infty.  It was a
startling consequence of the work of Casson and Gordon that in the
classical case the kernel of phi on C_1 is infinitely
generated. Beyond this, little has been discovered about the pair
(C_1,phi).

In this paper we present a new approach to studying C_1 by introducing
a group, M, defined as the quotient of the set of knots by the
equivalence relation generated by concordance and positive mutation,
with group operation induced by connected sum.  We prove there is a
factorization of phi, C_1-->M-->G_-.  Our main result is that both
maps have infinitely generated kernels.

Among geometric constructions on classical knots, the most subtle is
positive mutation.  Positive mutants are indistinguishable using
classical abelian knot invariants as well as by such modern invariants
as the Jones, Homfly or Kauffman polynomials.  Distinguishing positive
mutants up to concordance is a far more difficult problem; only one
example has been known until now.  The results in this paper provide,
among other results, the first infinite families of knots that are
distinct from their positive mutants, even up to concordance.}


\maketitlepage

\section{Introduction} The  classical knot concordance group,
$\calc$, was first defined in  \cite{fo}. Despite the remarkable
progress in low-dimensional topology since then, the structure of this
group remains largely a mystery. It is easily shown that
$\calc$ is countable and abelian. It is a consequence of Levine's work 
\cite{le1,le2} that there is a surjective homomorphism $\phi \!: \calc
\rightarrow \zz^\infty \oplus \zz_2^\infty \oplus
\zz_4^\infty$.  (Throughout this paper $G^\infty$ will denote an
infinite direct sum of groups and $\zz_n$ will denote the finite cyclic
group
$\zz/n\zz$.) Fox and Milnor 
\cite{fm} noted that
$\calc$ contains a subgroup isomorphic to
$\zz_2^\infty $; this subgroup is detected by
$\phi$. It follows from the work of Casson and Gordon  
\cite{cg1,cg2,ji}
that the kernel of $\phi$ contains a subgroup isomorphic to
$\zz^\infty$. Recently, Livingston  \cite{li1} described a subgroup of
the kernel isomorphic to
$\zz_2^\infty$. 

Levine proved that  the higher--dimensional analog of $\phi$ is
an isomorphism and hence provides a classification of
higher--dimensional knots up to concordance. Thus the homomorphism 
$\phi$ is the fundamental object in understanding  the
difference between classical knot concordance and the higher
dimensional theory.  

Beyond these results and their purely algebraic consequences, nothing
is known concerning the underlying algebraic structure of  $\calc$ or the
pair
$(\calc,\phi)$.  

Advances in geometric topology have offered new
insights into the geometric aspects of the knot concordance group. 
Freedman's work,
\cite{fr, fq}, shows that in the topological category all Alexander
polynomial one knots are trivial in
$\calc$.  The work of Donaldson \cite{do} and  Witten
\cite{witten} offer new obstructions in the smooth category to a knot
being trivial in $\calc$, \cite{cochrangompf, rudolph}.  Recently,
Cochran, Orr, and Teichner, \cite{cochranorrteichner}, have described a
nontrivial filtration of $\calc$, providing new obstructions to a knot
being trivial in $\calc$ in the topological category.  Yet, despite the
depth of these advances, the underlying structure of $\calc$ remains
mysterious.

Since the inception of the study of knot concordance, the role and
interplay of various 3--dimensional aspects of knotting and the
(basically 4--dimensional) study of concordance has been an ongoing
theme. The literature is extensive; topics include the
relationship of concordance to amphicheirality 
\cite{cm, go}, primeness  \cite{kl,li1}, invertibility 
\cite{li2}, hyperbolicity
 \cite{me}, mutation  \cite{ru}, and periodicity  \cite{na2, ck}. Of all
of these, positive mutation is, even from a purely 3--dimensional
perspective, the most subtle.

The difficulty in distinguishing a knot from its positive mutant 
 stems from the fact that positive mutation preserves the
$S$--equivalence class of a knot.  (This seems to be known, but a proof
has not appeared so one will be given in this paper.) As a
consequence, all abelian knot invariants---such as the Alexander
polynomial and module, signatures, and torsion invariants---are
incapable of distinguishing a knot from its positive mutant.
Furthermore, because the homology of the $n$--fold branched cyclic
cover of a knot and its positive mutant are isomorphic as
$\zz[\zz_n]$--modules, there is a correspondence between most of the
nonabelian representations that are used to distinguish knots when
abelian invariants fail. It is also known that  more recently
discovered knot invariants such as the Jones, Homfly, and Kauffman
polynomials also fail to distinguish a knot from its mutant, though
colored versions of these invariants are, at times, strong enough to
detect mutation \cite{mt}.  Other invariants that fail to distinguish a
knot from its mutant include its hyperbolic volume \cite{ru}, the eta
invariant  \cite{mr}, the hyperbolic torsion
\cite{po}, certain quantum invariants \cite{ro, li} and the
Casson-Walker invariants of surgery on the knot \cite{walker} and of
the 2--fold branched cover of the knot  \cite{mul}. 

The difficulties of studying the effect of mutation on concordance are
even deeper. In particular, since the algebraic concordance class of a
knot is determined by its Seifert form, it follows from the result that
positive mutants are $S$--equivalent that positive mutants are
algebraically concordant. That is, Levine's homomorphism cannot
distinguish a knot   from a positive mutant. Kearton  \cite{ker}
observed that the examples in \cite{li2}  demonstrated that negative
mutation can change concordance and soon after asked whether the same
was true for positive mutation \cite{k2}. The question was answered
affirmatively in \cite{kl2}, but as of yet the only result concerning
positive mutation and concordance is the single example presented in
\cite{kl2} of  two knots that are positive mutants but are not
concordant. This paper presents infinite classes of
such examples. 

In this article we provide a framework
for studying the interplay between concordance and  positive mutation
and then we identify some of the basic structures relating the two.

The results  are best described in terms of a new group,
the {\it Concordance Group of Mutants}, $\calm$, defined as the
quotient of $\calc$ obtained via the equivalence relation generated by
positive mutation. With the proof that positive mutants are
$S$--equivalent it will follow that Levine's homomorphism factors  as
$\phi_2 \circ \phi_1$, where $\phi_1 \!: \calc \rightarrow \calm$ and
$\phi_2 \!: \calm \rightarrow \zz^\infty \oplus
\zz_2^\infty
\oplus \zz_4^\infty$. Our main result is the following: 

\begin{thm} The kernels of $\phi_1$ and $\phi_2$ are infinitely
generated, containing subgroups isomorphic to
$\zz^\infty$.  The kernel of $\phi_2$ contains a subgroup isomorphic to
$\zz_2^\infty$.\end{thm} 

The argument is explicit; we construct  families of knots in
$S^3$ which generate appropriate subgroups of $\ker \phi_1$ and $\ker
\phi_2$.  The results are as follows. Let $K_a$ denote the 
$a(a+1)$--twisted double 
of the unknot. The case $a=2$ is drawn in Figure 2. 
\vskip3ex

\noindent{\bf Theorem \ref{phi2inf}}\qua   {\sl  Let $\{a_i\}_{i=1}^\infty$
be an increasing sequence of positive integers so that the sequence 
$\{2a_i+1\}_{i=1}^\infty$ is a sequence of primes. Then the knots
$K_{a_i}$ are linearly independent in
$\calm$ and generate a $\zz^\infty$ subgroup of   $\ker\phi_2\!:\calm\to
\zz^\infty
\oplus \zz_2^\infty \oplus \zz_4^\infty$.}
\vskip3ex

Figure 3 shows an order 2 knot $K_T$.  The ``$T$''  indicates that a
knot $T$ is to be tied in one band and $-T$ is to be tied in the other
band.  Let
$K_0$ denote the Figure 8 knot.

\vskip3ex

\noindent{\bf Corollary \ref{2torinf}}\qua  {\sl For the appropriate
choice of knots
$\{T_i\}_{i=1}^\infty$, the knots $K_{T_i}\# K_0$ generate an infinite
2-torsion subgroup of $\ker
\phi_2\!:\calm\to \zz^\infty
\oplus \zz_2^\infty \oplus \zz_4^\infty$.}\vskip3ex

Figure 4 shows a knot $K$ and two curves $B_1$, and $B_2^*$ in the
complement of $K$. Also shown is a 2-sphere intersecting $K$ in 4
points. Given a knot
$J$ let
$K_J^*$ denote the knot obtained by replacing the solid tori
neighborhoods of $B_1$ and $B_2^*$ by the exterior of $J$ and $-J$
respectively.  We give a recipe for choosing knots $J_i$ so that the
following theorem holds.
\vskip3ex

\noindent{\bf Theorem \ref{phi1ker}}\qua   {\sl There exists an infinite
collection of knots
$J_1, J_2,\ldots$ so that for any choice of integers $n_1,n_2,\ldots$
with only  finitely many of the $n_i$  nonzero, the connected sum
$ \#_i\ n_iK_{J_i}^*$  is not slice, but the knot obtained from 
$ \#_i\ n_iK_{J_i}^*$ by
performing a positive mutation on each summand 
$K_{J_i}^*$ using the indicated 2-sphere is slice. In particular the
kernel of
$\phi_1\!:{\cal C}\to {\cal M}$ contains a subgroup isomorphic to $
\zz^\infty$.
}

\vskip3ex

The main technique in proving these results is the use of Casson--Gordon
invariants. These are  sophisticated invariants which depend on an
auxiliary  vector in a  {\em metabolizer} for the linking form on
branched covers of the knot (a metabolizer is a  subspace  of the
homology of the branched cover on which the linking form vanishes). In
order to show a knot is not slice this method  requires that one
check that every metabolizer contains a vector for which some
appropriate Casson--Gordon invariant is non-trivial. Thus much of our
work involves a delicate examination of metabolizers. 

\vskip5ex

\noindent{\bf Outline}\qua This paper is divided into five sections. 
Following this introduction,  Section 2 presents background material.
In Section 2.1 we focus on 3--dimensional material, especially mutation
and the proof that positive mutation preserves the
$S$--equivalence class of a knot.  In fact, our proof shows that a knot
and its positive mutant admit Seifert surfaces with the same Seifert
form.  In Section 2.2 we review concordance and describe Levine's
homomorphism.  Casson--Gordon invariants are defined in Section 2.3 and
in Section 2.4 we discuss methods for their calculation. We use two types
of Casson--Gordon invariants, signatures and discriminants. Discriminants
are essential in working with
$\phi_1$ while signatures are needed to study $\phi_2$. We conjecture
that in general signature invariants cannot detect elements in the
kernel of $\phi_1$. 

Casson--Gordon invariants were first presented in two different ways: in
terms of intersection forms on infinite cyclic covers \cite{cg1} and in
terms of the limiting behavior of the signatures defined via
$n$--fold covers \cite{cg2}.  Here both approaches are needed. Since
the appearance of \cite{cg1,cg2} the approach of \cite{cg1} has become
standard and is almost always simpler and more conceptual. But for the
work of Section 3 we find that the approach of \cite{cg2} offers a much
clearer way to study the problem.

In Section 3 we prove that $\phi_2$ has an infinitely generated kernel,
containing a subgroup isomorphic
$\zz^\infty$. After summarizing the basic results of the section we
describe in Section 3.1 the needed algebra concerning metabolizers of
linking forms.  In the Section 3.2 we describe the geometry and
algebraic topology of the branched covers of mutant knots and in
Section 3.3 we give a new {\sl tangle addition} property for
Casson--Gordon invariants  (Theorem \ref{tangleadd}).  In much of the
work on Casson--Gordon invariants one uses
 Novikov additivity of the signature; to prove the tangle addition theorem
one must work in a context in which additivity fails, and we are forced
to look carefully at this failure of additivity via the Wall
nonadditivity result \cite{wa}.  In Section 3.4, we present the
necessary examples.  One step of the proof that these examples are
sufficient to show that the kernel of $\phi_2$  contains
$\zz^\infty$ is a lemma computing certain signatures of torus knots. 
This lemma is proved in Section 3.5.

In Section 4 we show that the kernel of $\phi_2$ also contains torsion,
in fact it contains a subgroup isomorphic to $\zz_2^\infty$.  That this
is the case for $\phi$ was shown in \cite{li3}.  We briefly review that
work here and demonstrate that the  infinite set of elements of
order 2 studied in \cite{li3} remain distinct when included in $\calm$.

Section 5 is devoted to proving that the kernel of $\phi_1$ contains a
subgroup isomorphic to $\zz^\infty$. Here a number of challenging
issues arise. The proof proceeds by starting with a family of slice
knots and showing that their positive mutants\eject are not slice. The proof
that the knots themselves are slice is nontrivial; it is carried out
in Section 5.1. The work of Sections 3 and 4 used signature invariants
associated to Casson--Gordon invariants.  In this section we use the
discriminant invariants instead.  Section 5.2 explores the covers of the
knots of interest.  In Section 5.3 we prove that the kernel of
$\phi_1$ is nontrivial and in the process develop the needed results
concerning the basic Casson--Gordon invariants associated with the
examples.  Sections 5.4 and  5.5 present the quite delicate analysis of
metabolizers of linking forms associated to the knots that is needed to
prove a linear independence result.  Section 5.6 is devoted to an
algebraic lemma.

The results of this paper generalize to the setting of knots in homology
spheres.  The proofs are essentially identical and the examples we give
apply in this more general setting. 

The authors wish to thank M Larsen for providing us with the proof of
Lemma
\ref{larsen} and J Davis for his help with representation theory.

\section{Background: Mutation, concordance and Casson--Gordon
invariants} 

\subsection{Classical knot theory and mutation} 

We will work in the smooth category, though, via \cite{fq}, all results
carry over to the topological locally flat setting. All homology groups
are taken with integer coefficients unless noted otherwise. 
For the basics of knot theory the reader is refered to \cite{rol}. 

In studying concordance, especially when working with such issues as
mutation, it is essential that orientation be carefully tracked. To do
so, one must begin with a precise definition: a knot consists of a
smooth oriented pair $(S,K)$ where
$S$ is diffeomorphic to the 3--sphere, $S^3$, and $K$ is diffeomorphic
to the 1--sphere, $S^1$. Knots are equivalent if they are oriented
diffeomorphic. With this, such notions as connected sum are defined as
they usually are for oriented manifolds and pairs.

This notation can become cumbersome and we will use it only when
necessary; we will usually follow the standard conventions,  writing
simply $K$ to denote a knot and viewing knots simply as oriented
1--spheres embedded in the standard (oriented) $S^3$. We will also
write $-K$ to denote the knot $(-S^3, -K)$, the ``mirror image'' of $K$
with orientation reversed. It is this knot that plays the role of
inverse in the concordance group.

\vskip3ex
\noindent {\bf Mutation}

Let $K$ be a knot in $S^3$ and let $S$ be a 2--sphere embedded in
$S^3$ meeting
$K$ transversely in exactly 4 points. Then $S$ bounds 3--balls,
$B_1$ and $B_2$, embedded in $S^3$, each of which intersects $K$ in two
arcs. The pairs $(B_i,B_i \cap K)$ are called tangles and we denote
them $T_1$ and $T_2$. Note that $K$ is recovered as the union $T_1
\cup_{id} T_2$, where $id$ denotes the identity map on
$S$. 

Suppose that $\tau$ is an orientation preserving involution of
$S$ with fixed point set (which necessarily consists of two points)
disjoint from $K \cap S$. Then we can form a new knot as the union $
K^* = T_1 \cup_{\tau} T_2$. This knot is called a mutant of
$K$. In the case that $\tau$ preserves the orientation of $S \cap K$ it
is called a positive mutant. 

(Though not used in this paper, it is worth noting here that given a 
knot $K$ and 2--sphere $S$ as above, the knot that results from
positive mutation is independent of the choice of $\tau$ in the sense
that any choice of orientation preserving involution which preserves
the orientation of $S\cap K$ gives the same result.
\cite{ru}.)

To each knot $K$ one can find a Seifert surface for it, say $F$, and to
$F$ along with a choice of basis for $H_1(F)$ one can associate a
Seifert matrix, say
$V$. Trotter \cite{tro2} defined a notion of
$S$--equivalence of Seifert matrices; roughly stated,
$S$--equivalence is an equivalence relation generated by change of
basis and stabilization of a certain type. A basic results states that
any two Seifert matrices for a given knot are
$S$--equivalent. Knots are called $S$--equivalent if they have
$S$--equivalent Seifert matrices. We now outline a proof of a result
that implies that positive mutation preserves the
$S$--equivalence class of a knot.

\begin{thm} If a knot $K^*$ is a positive mutant of $K$, then there are
Seifert surfaces $F$ and $F^*$ for $K$ and $K^*$, respectively, so
that, for an appropriate choice of bases, the Seifert matrices are the
same.  In particular, $K$, and $K^*$ are $S$--equivalent. \end{thm}

\begin{proof} To set up notation, let $K$ intersect $S$ in points $p_+$,
$p_-$, $q_+$, and $q_-$, where $p_+$ is joined to $p_-$ and $q_+$ is
joined to
$q_-$ via oriented arcs $A^1_1$ and $A^1_2$ in $T_1 = K \cap B_1$.
Similarly, $p_+$ is joined to $q_-$ and $q_+$ is joined to
$p_-$via arcs $A^2_1$ and $A^2_2$ in $T_2 = K \cap B_2$.

Fix a pair of arcs on $S$, $A^3_1$ and $A^3_2$, joining $p_+$ to
$p_-$ and
$q_+$ to $q_-$ on $S$ such that $\tau (A^3_1) = A^3_2$. Then the pair of
oriented circles $A^1_1 \cup A^3_1$ and $A^1_2 \cup A^3_2$ is the
oriented boundary of a properly embedded oriented connected surface
$F_1$ in $B_1$. The proof is a standard argument in relative
obstruction theory, considering maps of the complement of the arcs
($A_1^1$ and $A_2^1$) to $S^1$ and letting $F_1$ be the preimage of a
regular value. Similarly, the union of arcs
$A^2_1 \cup A^3_2 \cup  A^2_2 \cup A^3_1$ forms a circle, and this
circle is the oriented boundary of a properly embedded oriented
connected surface
$F_2$ in
$B_2$. The union $ F = F_1 \cup F_2$ is an oriented Seifert surface for
$K$. Since $\tau$ is orientation preserving and
$\tau(A^3_1) = A^3_2$, $F_1$ and $F_2$ can similarly be joined to form
an oriented Seifert surface $F^*$ for $K^*$. 

Pick oriented bases for $H_1(F_1)$ and $H_1(F_2)$. A basis for
$H_1(F)$ can be formed as the union of these, along with one more
element that can be represented by a curve on $F$ that intersects
$A^3_1$ and $A^3_2$ each exactly once. 

It is now clear that there will be a corresponding basis for
$H_1(F^*)$. Furthermore, working directly from the definition of the
Seifert matrix, one sees that the corresponding Seifert matrices are
identical as follows. The pairings between elements of the basis
represented by curves that miss $S$ are clearly unchanged. Let $\alpha$
be a curve representing the extra element of $H_1(F)$ that intersects
$A^3_1$ and $A^3_2$ each exactly once. It is clear that the linking
number of $\alpha$ with the positive push-off of any curve missing $S$
is unchanged (here the fact that the mutation is positive is
essential). The self-linking of
$\alpha$ with its positive push-off is unchanged (regardless of the
sign of the mutation). \end{proof}

\subsection{Concordance, mutation, and the Levine homomorphism}

\noindent{\bf Concordance}

Returning to our more formal notation momentarily, a knot $(S,K)$ is
called slice if there is a proper 4--manifold pair $(B,D)$ with
$B$ diffeomorphic to the 4--ball, $B^4$, $D$ diffeomorphic to the
2--disk, $B^2$, and with $\partial(B,D) = (S,K)$. Knots $K_1$ and
$K_2$ are called concordant if the connected sum $K_1 \# -K_2$ is slice.
Concordance is an equivalence relation and the set of equivalence
classes forms a group under the operation induced by connected sum.
This group is called the concordance group, denoted $\calc$. 

An important alternative approach to concordance is the following.
Knots $(S_1,K_1)$ and $(S_2,K_2)$ are called concordant if there is a
properly embedded annulus,  $A$ (a diffeomorph of $S^1 \times [0,1]$), in
$S^3 \times [0,1]$ with
$\partial(S^3 \times [0,1], A) = (S_1,K_1) \cup -(S_2,K_2)$. In this
case $A$ is called a concordance from $(S_1,K_1)$ to
$(S_2,K_2)$. 

\vskip2ex
\noindent{\bf The algebraic concordance group and Levine's
homomorphism} 

The definition of Levine's algebraic concordance group \cite{le1} in
dimension 3 proceeds as follows. Consider the set of all even
dimensional integer matrices
$V$ satisfying $\det(V - V^t ) =1$. Such a matrix is called metabolic
if it is congruent to a matrix having a half dimensional block of
zeroes in its upper left hand corner. Matrices $A_1$ and $A_2$ are
called algebraically concordant if the block sum, $A_1 \oplus -A_2$ is
metabolic. In particular
$S$--equivalent matrices are algebraically concordant.  This defines an
equivalence relation and the set of equivalence classes forms a group
under block sum called the algebraic concordance group, denoted
$\calg$. Levine proves in \cite{le2} that $\calg \cong
\zz^\infty \oplus
\zz_2^\infty \oplus
\zz_4^\infty $.

There is a homomorphism from $\calc$ to $\calg$, denoted $\phi$,
induced by the map that assigns to a knot an associated Seifert matrix.
Levine proves that this is surjective. Its analogue in
higher dimensions is injective and either surjective or onto an index 2
subgroup. In dimension 3 it is known that the map is onto, and Casson
and Gordon proved that it is not injective. In fact, by \cite{ji} and
\cite{li3} the kernel of Levine's homomorphism contains a subgroup
isomorphic to
$\zz^\infty  \oplus
\zz_2^\infty $.

\vskip3ex
\noindent{\bf The mutation group}

The two relations, concordance and positive mutation, generate an
equivalence relation on the set of knots. Denote the set of equivalence
classes by
$\calm$. Thus the two knots, $K_1$ and $K_2$, represent the same element
of $\calm$ if a finite sequence of concordances and positive
mutations converts one into the other.

It is clear that if $K_1$ and $K_2$ represent the same class in
$\calm$ and $J_1$ and 
$J_2$ represent the same class in $\calm$,  then $K_1 \# J_1$ and
$K_2
\# J_2$ represent the same class in $\calm$. Hence, connected sum
induces an operation on
$\calm$. Clearly, the unknot serves as an identity under this operation
and mirror images provide inverses, so $\calm$ is a group. It is easily
seen that this is quotient of $\calc$. Finally, since by Theorem 2.2
mutation does not alter the image of a knot under Levine's
homomorphism, we have the following result. 

\begin{thm} Levine's homomorphism $\phi\!:\calc \rightarrow
\calg$ factors as the composition of $\phi_1\!:\calc \rightarrow
\calm$ and $\phi_2\!:\calm
\rightarrow \calg$.\end{thm}

\eject
\subsection{Casson--Gordon invariants}\label{cgsection}
\noindent{\bf Casson--Gordon invariants as obstructions to slicing} 

We now present the formulations of Casson--Gordon invariants that will
be needed. Let
$K$ be a knot in $S^3$ and let $M_d$ denote its $d$--fold cyclic
branched cover. In what follows $d$ will always be a prime power,
$q^n$, and in particular,
$H_1(M_d)$ will always be finite; in fact, $M_d$ is a
$\zz_q$--homology sphere \cite{cg1}. There is a symmetric
  {\em linking form}  on the finite abelian group
$H_1(M_d;\zz)$,  
$$lk\!: H_{1}(M_d;\zz)
\times H_{1}(M_d;\zz) \rightarrow \qq/\zz.$$ This linking form is
nonsingular in the sense that it defines an isomorphism
$H_1(M_d;\zz) \rightarrow$ Hom$(H_1(M_d;\zz),\qq /\zz)$.

\begin{defn} A {\em metabolizer}
for the linking form on $H_1(M_d;\zz)$ is a subgroup $A
\subset H_1(M_d;\zz)$ for which $A = A^{\perp}$. By nonsingularity we
have
$|A|^{2} = |H_1(M_d;\zz)|$.
\end{defn}

Next, let $p$ be prime power, and let $\chi$ be a character,
$\chi
\!: H_1(M_d)
\rightarrow \zz_p$. We will define two invariants in this setting,
$\delta(K,\chi)$ and $\sigma_b(K,\chi,i)$. (The subscript  $b$ stands  
for bounded, and $i$ is an integer, $0<i<p$. The invariants $\delta$ and
$\sigma_b$  can be  defined for arbitrary prime power covers of
$K$, but we will only define and work with $\sigma_b$ for the
$2$--fold covers, and with $\delta$ for the 3--fold covers.) The main
theorem, derived from the work of Casson and Gordon in
\cite{cg1,cg2} is the  following.  

\begin{thm}\label{cassongordon} If $K$ is slice there is a metabolizer
$A$ of
$H_1(M_d)$  which is invariant under the action of the group of deck
transformations acting on $H_1(M_d)$ so that, for any $\chi \!:
H_1(M_d) \to \zz_p$ that vanishes on $A$, one has
$\sigma_b(M_d,\chi) =0$ and
$\delta(M_d,\chi) =1$. \end{thm} 

Note that $\sigma_b(M_d,\chi)$ will be seen to take values in an
additive group and
$\delta(M_d,\chi) $ will take values in a multiplicative group, and
that is why one takes value 0 and the other takes value 1. The proof of
this theorem is contained in the rest of this subsection after presenting
the relevant definitions.

\vskip2ex
\noindent {\bf Definition of the Casson--Gordon invariant
$\sigma_b$ }

To begin we need to review bordism theory.  A good reference for details
is \cite{cf}. 
The $n$-dimensional bordism group of a group $G$, denoted $\Omega_n(G)$,
is defined to be the set of equivalence classes of pairs $(M^n, \chi)$
where $M^n$ is a connected oriented closed manifold and $\chi$ is a
homomorphism from $\pi_1(M^n)$ to $G$; two such pairs, $(M_1^n, \chi_1)$
and $(M_2^n,
\chi_2)$ are considered equivalent if there is an oriented and compact
$(n+1)$ manifold
$W^{n+1}$ with boundary the disjoint union of $M_1^n$ and $-M_2^n$ and a
homomorphism from $\pi_1(W^{n+1})$ to $G$ restricting to give $\chi_1$ and
$\chi_2$ on the boundary.  The group structure is induced by the connected
sum operation; the group is abelian.

Basic results in bordism theory give that for the trivial group $0$,
$\Omega_0(0) \cong \zz$, and $\Omega_1(0) \cong \Omega_2(0)
\cong \Omega_3(0) \cong
0$.  There is a bordism spectral sequence with $E^2$ term given by $H_i(G,
\Omega_j(0))$. It follows readily that for $n \le 3$, the natural map
$\Omega_n(G) \rightarrow H_n(G)$ is an isomorphism.  In particular,
$\Omega_3(\zz_p) \cong \zz_p$.

To apply this in our situation, let $\chi_1$ be a $\zz_p$--valued
character defined on
$H_1(M_2)$. For each   integer $k$ there is a natural projection map
$H_1(M_{2^k}) \rightarrow H_1(M_2)$, and hence   an induced character
$\chi_k \! :H_1(M_{2^k}) \rightarrow \zz_p$. 
Thus $p$ copies of $(M_{2^k},\chi_k)$ bound a pair $(W_{k},\chi_k)$.
More precisely, since $\Omega_3(\zz_p)\cong \zz_p$, $p(M_{2^k},\chi_k)=0$
in
$\Omega_3(\zz_p)$. This means that there is a 4--manifold and character
$(W_{k},\chi_k')$ whose boundary is $p$ disjoint copies of $M_{2^k}$ and
with
$\chi_k'$ restricting to $\chi_k$ on each boundary component. Let
$\tilde{W}_{k}$ denote the corresponding
$p$--fold cover, with group of deck transformations isomorphic to
$\zz_p$, generated by
$R_k$.

Under the action of $R_k$, the homology group
$H_2(\tilde{W}_k,\cc)$ splits into eigenspaces
$H_{k,i} = H_2(\tilde{W}_k,\cc)_{\zeta^i}$, where $\zeta$ is a fixed
primitive
$p$th root of unity, $1 \le i \le p-1$. Let $\sigma(K,\chi,k,i) =
\frac{1}{p} (\mbox{signature}(H_{k,i}) -
\mbox{signature}(H_2(W_k,\cc)))$, where the signatures are those of the
intersection forms of the 4--manifolds, properly restricted. 

At this point one can consider the sequence $\sigma_b(K,\chi,i) =
\{ \sigma(K,
\chi, k, i) \}_{k=1, \ldots ,\infty}$. Casson and Gordon prove in
\cite{cg1}: 

\begin{thm}\label{cg2} If $K$ is slice, there is a metabolizer $A$ of
$H_1(M_2;\zz)$ such that if $\chi$ is a $\zz_p$--valued character
vanishing on
$A$, then for each $i$, the sequence $\{ \sigma(K,\chi,k,i) \}_k$ is
bounded.
\end{thm}

As a result, if we view the sequence as in
$(\prod_{1}^\infty\qq)_b$, the set of all infinite sequences modulo
bounded sequences, we have
$\sigma_b(K,\chi,i) = 0\in (\prod_{1}^\infty\qq)_b$ and the
$\sigma_b$ part of Theorem \ref{cassongordon} follows. 

\vskip2ex

\noindent {\bf Definition of the Casson--Gordon invariant $\delta$ } 

The Casson--Gordon invariant $\delta$ was first defined via a
discriminant of a Casson--Gordon Witt class valued invariant
$\tau$ defined in \cite{cg1}. See \cite{lit1,gl2} for details. However,
it was shown in \cite{kl1} that if one is interested only in the
discriminant, the 4--dimensional work of defining $\delta$ can be
bypassed and there is a simple 3--dimensional interpretation in terms
of twisted Alexander polynomials.

Let $\bar{M}_d$ denote the 3--manifold obtained from $M_d$ by removing
the lift of $K$ in the branched cover; that is,
$\bar{M}_d$ is the $d$--fold cyclic cover of the complement of
$K$. There is an inclusion of $H_1(\bar{M}_d)
\rightarrow H_1({M}_d) $ and hence an induced character
$\chi'\!:H_1(\bar{M}_d)
\rightarrow
\zz_p$. Let $\tilde{M}$ denote the infinite cyclic cover of the
complement of
$K$. It maps to $\bar{M}_d$ and hence there is a character
$\tilde{\chi}\!:H_1(\tilde{M})
\rightarrow
\zz_p$.

Suppose that $d$ and $p$ are odd prime powers.  Letting $\zeta_p$ denote
a primitive
$p$th root of unity, one can consider the twisted homology group
$H_1(\tilde{M},\qq(\zeta_p))$ where the twisting is via
$\tilde{\chi}$. Now $\zz$ acts on $\tilde{M}$ via deck transformations,
where
$\tilde{M}$ is viewed as the infinite cyclic cover of $M_d$. Hence the
homology group $H_1(\tilde{M},\qq(\zeta_p))$ is a
$\qq(\zeta_p)[t,t^{-1}]$ module. This ring is a PID. For the moment let
$D$ denote the order of that module (which is a torsion module
according to Casson and Gordon \cite{cg1}). This is an element of
$\qq(\zeta_p)[t,t^{-1}]$, well defined up to units. The discriminant
Casson--Gordon invariant
$\delta(M,\chi)$ is defined to be the image of $D$ in the the
multiplicative group
$\qq(\zeta_p)[t,t^{-1}]^{\times}/N$, where the $\times$ denotes nonzero
elements and $N$ denotes the subgroup consisting of norms: elements of
the form $f(t)\overline{f(t)}$ where the bar acts by complex
conjugation on the coefficients and by sending $t$ to
$t^{-1}$. 

The theorem proved in \cite{kl1} from which the second half of Theorem
\ref{cassongordon} follows is the following.

\begin{thm} If $K$ is slice, there is a $\zz_d$--invariant metabolizer
$A$ for
$H_1(M_d)$ such that for all characters $\chi$ vanishing on $A$,
$\delta(M_d,\chi) = 1 \in \qq(\zeta_p)[t,t^{-1}]^{\times}/N$.
\end{thm} 

\subsection{Computation of Casson--Gordon invariants} 

\noindent{\bf Additivity of Casson--Gordon invariants} 

Given knots $K_1$ and $K_2$  and characters $\chi_1$ and $\chi_2$ on
the appropriate
$d$--fold covers, one has an additivity property for Casson--Gordon
invariants. Such a result was first proved in \cite{ji,gi}. In the
present context the result follows almost immediately from the
definition. 

\begin{thm} \label{sigmai}In the setting just described, $$\sigma_b(K_1
\# K_2,\chi_1 \#
\chi_2) = \sigma_b(K_1 ,\chi_1 ) + \sigma_b( K_2, \chi_2)$$ and
$$\delta(K_1 \# K_2,\chi_1 \#
\chi_2) =  \delta(K_1 ,\chi_1 ) \cdot \delta( K_2, \chi_2).$$\end{thm}

\noindent{\bf The geometry of satellite knots and their covers} 

In general the computation of Casson--Gordon invariants for a particular
knot can be quite difficult. However, it is often sufficient to
understand how the invariants change under particular geometric
operations. Roughly stated, if a geometric construction on a knot
yields a slice knot, then it is known that certain Casson--Gordon
invariants vanish, so the invariants of the original knot are
determined by the effect of the geometric construction. This underlying
theme has been central to most calculations since it was first used in
\cite{gi,lit2}.

The construction that we are most interested in is that of forming a
satellite knot from a given knot. The following perspective on this
operation was first described in \cite{lit2}.

Let $U$ be an unknot in
$S^3$ and let
$J$ be a knot in
$S^3$. Removing a tubular neighborhood of $U$ and sewing in the
complement of
$J$ results in $S^3$ again, provided that the meridian and longitude of
$J$ are mapped to the longitude and meridian of $U$, respectively. (The
resulting space is the union of a knot complement and a solid torus
that ``fills in'' the knot.)

If $K$ is a knot in the complement of $U$, then after this operation
$K$ represents a perhaps different knot, $K'$, in $S^3$. In effect,
viewed in the complement of $U$,
$K$ represents a knot in a solid torus. The knot $K'$ is formed by
tying that solid torus into the knot $J$.  (The standard terminology
is to call $K'$ a {\it satellite knot} with {\it companion} $J$ and
{\it embellishment} $K$.) 
\vskip2ex
\noindent{\bf The effect of the satellite operation on
$\sigma_b$}

For a  different form of Casson--Gordon invariants the connection
between the invariants and satellite knots was worked out by Litherland
\cite{lit2}. Hence we only present an outline of the approach and give
the details as they differ from \cite{lit2} and are needed here. 

We will suppose that $U$ is nullhomologous in the complement of
$K$. Hence, the 2--fold branched cover of $K'$, $M_2'$, is obtained
from the the 2--fold branched cover of
$K$ by removing the two lifts of $U$ from $M_2$ and replacing each with
the complement of $J$. This has no effect on the homology, and hence
there is a correspondence between
$\zz_p$--valued characters $\chi$ on $H_1(M_2)$ and $\zz_p$--valued
characters on
$H_1(M_2')$.

Suppose that $\chi$ takes value $a$ on one lift of $U$. Then it takes
value $-a$ on the other lift, since the deck transformation of the
2--fold cover of a knot acts by multiplication by $-1$ on the first
homology. 

Recall that the (Levine--Tristram) signature function $\si_t(K),
t\in \rr$ of a knot
$K$ is the signature of the symmetric matrix 
$$B(e^{2\pi i t})=(1-e^{2\pi i t})V + (1-e^{-2\pi i
t})V^T\eqno{(1)}$$ where $V$ is a Seifert matrix for $K$. Note that 
$\si_t(K)=\si_{-t}(K)=\si_{t+n}(K)$ for $n\in\zz$, so that
$\si_t(K)$ is determined by its values for
$t\in[0,{1\over 2}]$ and also that $\si_t(K)=-\si_t(-K)$.

Now, considering the $2^k$--fold covers, $M_{2^k}'$ is obtained from
$M_{2^k}$ by removing the $2^k$ lifts of $U$ and replacing them with
copies of the complement of
$J$. Under the character on these covers corresponding to $\chi$, half
of the lifts of $U$ map to $a$ and half map to $-a$. This motivates the
following result. Notice that in the statement of this theorem, the
term $(2^k \sigma_{ja /p}(J))$ represents a sequence in
$(\prod_{1}^\infty\qq)_b$, the group of infinite sequences modulo
bounded sequences. 

\begin{thm}\label{lith} If $U$ is unknotted and null homologous in the
complement of
$K$ and $\chi$ and $\chi'$ are corresponding characters on the covers
$M_2$ and
$M_2'$ just described, then $\sigma_b(K,\chi,j) -
\sigma_b(K',\chi',j)  = 2^k
\sigma_{ja / p}(J)$.\end{thm}

\begin{proof} The idea of the proof is simple and follows the approach of
\cite{lit2} or \cite{gl1}. It is actually implicit in \cite{gi}.
Basically, if $W$ denotes a 4--manifold used in the computation of
$\sigma_b(K,\chi,j)$ then a 4--manifold that can be used to compute
$\sigma_b(K',\chi',j)$ is built from $W$ by adding copies of a
4--manifold $Y$ along the lifts of $U$. The manifold $Y$ is a
4--manifold with boundary 0--surgery on $J$, $N$, for which the natural
map of $H_1(N) \rightarrow \zz$ extends. But it is this manifold $Y$
that can be used to compute the classical signatures of
$J$. Notice that since $\sigma_{ja / p}(J) = \sigma_{-ja / p}(J)$ the
sign issue disappears. \end{proof}

\noindent{\bf The effect of the satellite operation on $\delta$} 

In the examples we will be considering in which $\delta$ must be
evaluated the situation will be somewhat different. In these cases $U$
will not be null-homologous in $S^3 - K$; in fact, $U$ will represent
$d$ times a generator of
$H_1(S^3 - K)$ and we will be considering the $d$--fold cover of the
knot. (We can fix an orientation of $U$ so that $d$ is positive, and
use this orientation to orient the lifts of $U$ in what follows.) In
this case, $U$ lifts to
$d$ curves, say
$\{U_i\}_{i = 1,
\ldots, d}$, where the deck transformation $T$ cyclicly permutes the
$U_i$. 

If we consider a $\zz_p$--valued character $\chi$ on $H_1(M_d)$ and
$\chi$ takes value $a$ on $U_1$, then the value of $\chi$ on
$U_i$ will be ${T^*}^i(\chi)(a)$, where $T^*$ is the transformation on
$H^1(M_d;\zz_p)$ induced by $T$. To simplify notation we define $a_i =
{T^*}^i(\chi)(a)$. 

In defining $\delta$ we also needed to consider the representation of
$H_1(\bar{M_d})$ to $\zz$. (Recall that $\bar{M_d}$ is the cyclic cover
of $S^3 - K$.) In the setting we just describe, it is clear that this
representation takes value 1 on each of the
$U_i$. With this we have: 

\begin{thm}\label{deltathm} In the setting just described, $\delta(K') /
\delta(K) = \prod_{i=1}^d \Delta_J(\zeta^{a_i}t)$, where
$\Delta_J(t)$ is the Alexander polynomial of $J$. \end{thm}

\begin{proof} The idea of the proof is related to that of the previous
theorem. In the definition of $\delta$ one considers the infinite
cyclic cover of
$\bar{M_d}$. The infinite cyclic cover of $\bar{M_d'}$ is built from
that of
$\bar{M_d}$ by removing $d$ copies of $\rr \times B^2$ and replacing it
with $d$--copies of the infinite cyclic cover of $J$. Hence it is
expected that $d$ factors of the Alexander polynomial of $J$ should
appear. The $\zeta^{a_i}$ appear because of the twisting in the
homology. That the product appears as it does follows from a
Mayer-Vietoris type argument. Details are presented in
\cite{kl}.\end{proof}

\section{The kernel of $\phi_2$ contains $\zz^\infty$} 

We begin with a definition.

\begin{defn} A knot $K$ is called {\em cg--slice} if there exists a
metabolizer $A$ for $H_1(M_2)$ such that for all prime order characters
$\chi$ on
$H_1(M_2)$ that vanish on $A$, and for all i, $\sigma_b(K,\chi,i) =
0$.\end{defn}

Hence, the following is a restatement of Theorem \ref{cg2}:

\begin{thm}\label{abc} If $K$ is slice, then it is cg--slice. \end{thm}

Therefore it follows that

\begin{cor} If $J$ is concordant to $K$, then
$K
\# -J$ is cg--slice.\end{cor}

We first prove, in Section 3.1, the following: 

\begin{thm}\label{addcg} If $K_1$ is cg--slice and $K_1 \# K_2$ is
cg--slice, then $K_2$ is cg--slice.\end{thm}

The proof is basically algebraic, based on a careful study of
metabolizers. The argument is similar to that given by Kervaire
\cite{ke} in the vector space (as opposed to finite group) setting.

The key result of this section is proved in Section 3.3:

\begin{thm}\label{posmut} If $J$ is a positive mutant of $K$ then $K
\# -J$ is cg--slice.\end{thm}

This result follows from  a technical result, a ``tangle addition''
property for $\sigma_b$.

As an immediate consequence of these four results we have the corollary
that is used to analyze the kernel of
$\phi_2$:

\begin{cor}\label{posmutcor} If $K$ represents the trivial class in the
concordance group of mutants,
$\calm$, then $K$ is cg--slice.\end{cor}

Section \ref{examps1} considers examples. We note first that the
original examples in \cite{cg2} can be used to produce the desired
elements in
$\mbox{Ker}(\phi_2)$. We then explore these examples further. 

\subsection{Cancellation and metabolizers} The proof of \ref{posmut}
depends an algebraic result, modeled on one presented by Kervaire
\cite{ke} in the vector space setting.

Let $G_1$ and $G_2$ be finite abelian groups with linking forms
(nonsingular symmetric
$\qq /
\zz$--valued forms), denoted with brackets: $\langle\ ,\ \rangle\!: G_i
\times G_i
\rightarrow
\qq/\zz$. Suppose that
$A_1$ is a metabolizer for $G_1$ and $A$ is a metabolizer for
$G_1 \times G_2$ with respect to the product linking form. Define
$A_2$ by: $$A_2 = \{g \in G_2 | (g_1,g) \in A \mbox{ for some } g_1 \in
A_1\}.$$ 

\begin{thm} $A_2$ is a metabolizer for $G_2$.\end{thm} 

\begin{proof} First note that  a subgroup $B$ of $H$, a torsion group with
nonsingular linking pairing, is a metabolizer if and only if the
linking form vanishes on
$B$ and the order of $B$ is the square root of the order of $H$. (From
the exact sequence $0 \rightarrow B^\perp
\rightarrow H \rightarrow \mbox{Hom}(B,\qq/\zz) \rightarrow 0$ we have
that $|B||B^\perp| = |H|$. If the linking form vanishes on
$B$, then $B \subset B^\perp$.) 

It is clear that the linking form for
$G_2$ vanishes on
$A_2$. All that we need to show is that the order of $G_2$ is the square
of the order of
$A_2$. 

To simplify notation we will view $G_{1}$ and $G_{2}$ as subgroups of
$G_{1}\times G_{2}$ in the natural way. Let $p$ denote the projection
from $G_{1}\times G_{2}$ to
$G_{1}$. Finally, define $A_{0} = \{(x,y) \in A | x \in A_{1}\}$. 

We have the two exact sequences:
$$0 \rightarrow A_{0} \rightarrow A \rightarrow p(A)/(A_{0} \cap p(A))
\rightarrow 0$$
$$0 \rightarrow A_{1} \cap A \rightarrow A_{0} \rightarrow A_{2}
\rightarrow 0.$$
Notice that $A_{1} \cap A$ pairs trivially against both $p(A)$ and
$A_{1}$. Hence it pairs trivially against the sum, $p(A)+ A_{1}$. In
particular, since the pairing is nonsingular, we have that the orders
satisfy
$|A_{1} \cap A|\ |p(A)+ A_{1}| \le |G_{1}|$. Notice also that the order
of the sum of groups is equal to the product of their orders, divided by
the order of the intersection. Hence 
$${|A_{1} \cap A|\ |p(A) |\ |A_{1}| \over \ |p(A) \cap A_{1}|}
\le |G_{1}|.$$
From the second exact sequence we have $$|A_{2}| = {|A_{0}|
\over |A_{1} 
\cap A|}.$$ From the first exact sequence it follows that
$$|A_{2}| = {|A||A_{1} \cap p(A)| \over |A_{0} \cap A| |p(A)|}.$$
Using the inequality (along with the fact that $A_0 \cap p(A) = A_1
\cap p(A)$) applied to the denominator, we have $$|A_{2}|
\ge {|A| |A_{1}| \over |G_{1}|}.$$ Since each metabolizer has order the
square root of the order of its ambient group, and since no 
self--annihilating subgroup can have order greater than the square
root of the order of the ambient group, the result follows. \end{proof}

We can now prove Theorem \ref{addcg}.

\begin{proof}  If
$K_1$ is cg--slice, let
$A_1$ be the appropriate metabolizer. Similarly, let $A$ be the
metabolizer that exists because
$K_1 \# K_2$ is cg--slice. Then the previous theorem provides a
metabolizer $A_2$. For any $\chi_2 \in A_2$ there is a $\chi_1
\in A_1$ such that $(\chi_1,\chi_2) \in A$. Hence, $0 =
\sigma_b(K_1 \# K_2) = \sigma_b(K_1 ) +
\sigma_b( K_2) =
\sigma_b( K_2)$. Hence, $K_2$ is cg--slice as desired. \end{proof}

\subsection{Decomposing the  branched covers of mutant knots}

\vskip.2in

\noindent{\bf The regular cyclic branched cover} 

In order to prove a
tangle addition theorem for Casson--Gordon invariants in the next
section, we must first understand the structure of the relevant covers
of a knot and its positive mutant.  In particular, a detailed
understanding of the gluing maps used to assemble the cover as the
union of covers of tangles is needed.  We develop the necessary results
in this section.

If $K$ and $K^*$ are positive mutants, let $S$ be the splitting
2--sphere along which the mutation is formed, yielding the tangle
decompositions 
$K= T_a
\cup_{id} T_b$ and $K^*= T_a
\cup_{\tau} T_b$, where $\tau$ is an involution of $S$.

The $2^k$--fold cyclic cover of the complement of $K$ corresponds to
the kernel of the natural map of $\pi_1(S^3 - K) \rightarrow
\zz_{2^k}$. Since $\tau$ preserves the kernel of the corresponding map
$\pi_1(S - \{S \cap K\}) \rightarrow
\zz_{2^k}$, and has a fixed point, by elementary covering space theory
it lifts to an involution with fixed point on the $2^k$--fold branched
cover of $S$, denoted $F_k$.

\begin{thm}  $\tau_k$ acts by multiplication by -1 on
$H_1(F_k)$.\end{thm}

\begin{proof} The quotient of $F_k$ under the $\zz_2$--action of
$\tau_k$ is a
$2^k$--fold branched cover of the quotient of $S$ under the action of
$\tau$.  This quotient is a 2--sphere and the branch set has quotient
consisting of  2 points.  The
$2^k$--fold branched cover of a $S^2$ over 2 points is again a
2--sphere, so the quotient of $F_k$ under $\tau_k$ is a 2--sphere. 

Since $\tau_k$ is an involution, $H_1(F_k,\rr)$ splits as a direct sum
of the
$+1$ and
$-1$ eigenspaces.  A transfer argument (see \cite{bredon}) shows that the
+1 eigenspace is isomorphic to the first homology of the quotient, which
is trivial since the quotient is a sphere.  The theorem follows. \end{proof}

\noindent{\bf The irregular branched cover}

Since the Casson--Gordon invariants we are interested in are associated
with the $p$--fold cyclic covers of the $2^k$--fold branched cyclic
covers of
$K$ and $K^*$, we must next consider the decomposition of these
spaces.  

The $p$--fold cover of $M_k$ induced by $\chi$ decomposes as
$\tilde{M_k^a}\cup_{id}
\tilde{M_k^b}$, where the union is along $\tilde{F_k}$, the induced
$p$--fold cover of $F_k$.  Similarly, the $p$--fold cover of
$M_k^*$ induced by $\chi_k$ decomposes as  $\tilde{M_k^a} 
\cup_{\tilde{\tau}_k} 
\tilde{M_k^b}$ where $\tilde{\tau}_k$ is some lift of $\tau_k$. Notice
that $\tilde{\tau}_k$ leaves invariant the preimage of a fixed point of
$\tau$ and, by composing with a deck transformation of
$\tilde{M_k^b}$  if need be, we can assume that $\tilde{\tau}_k$ fixes
a  point in that set.

It might be the case that
$\tilde{F_k}$ is disconnected and this makes the covering space
arguments a bit more delicate, since one loses some of the uniqueness
properties of liftings that exist in the connected case. In particular,
for connected covers one has that if a deck transformation fixes a
point it is the identity; this is not longer true for disconnected
covers.  

\begin{thm}  If $R_k$ is the restriction of the order $p$ deck
transformation of $\tilde{M_k}$ to $\tilde{F_k}$, then
$\tilde{\tau}_k$ and $R_k$ together generate a dihedral action on
$\tilde{F_k}$. \end{thm}

\begin{proof} The cover of $M_k^*$ decomposes as  $\tilde{M_k^a} 
\cup_{\tilde{\tau}_k} \tilde{M_k^b}$, but note that the cyclic action
of the deck transformation restricted to  $\tilde{M_k^b}$ is the
inverse of the action that arises when we consider the cover of $M_k$,
since the character $\chi$ has been replaced by $-\chi$.  That the
$\zz_p$ action is well defined on the union implies that
$\tilde{\tau}_k \circ R_k = R_k^{-1} \circ \tilde{\tau}_k$.  

Consider the restriction of $\tilde{\tau}_k$ and $R_k$ to the preimage
of a fixed point of $\tau_k$, a set with $p$ elements.  Then $R_k$ acts
as a $p$--cycle and $\tilde{\tau}_k^2$ commutes with $R_k$.  In the
symmetric group on $p$ letters the only elements that commute with a
$p$--cycle $\rho$ are powers of $\rho$.  But
$\tilde{\tau}_k^2$ has a fixed point  and the only power of
$\rho$ with a fixed point is the identity.

We now have the relations $\tilde{\tau}_k^2 = 1, R_k^p = 1$, and
$\tilde{\tau}_k \circ R_k = R_k^{-1} \circ \tilde{\tau}_k$.  Since
neither 
$\tilde{\tau}_k$ or $R_k$ are trivial, the result follows.
\end{proof}

\noindent{\bf The action on homology}  

We now want to understand the action of $\tilde{\tau}_k$ and $R_k$ on
the homology of $\tilde{F_k}$.  At this point we work with
$H_1(\tilde{F}_k,\cc) \cong H_1(\tilde{F}_k) \otimes \cc$.

First note that since $R_k$ satisfies $R_k^p = id$, we have that 
$H_1(\tilde{F}_k,\cc)$ splits into $\zeta^i$ eigenspaces ($\zeta$ is a
fixed primitive $p$ root of unity).  We denote these eigenspaces
$H_{\zeta^i}, i = 0, \ldots , p-1$.  

Since $H_1(\tilde{F}_k,\cc) \cong H_1(\tilde{F}_k)
\otimes \cc$,  complex conjugation induces an action on
$H_1(\tilde{F}_k,\cc)$.  Via averaging over the action of
$R_k$, we see that every class in $H_{\zeta^i}$ can be represented as a
sum of classes of the form $\sum_{j=0}^{p-1} R_k^j x \otimes \alpha
\zeta^{-ij}$.  It follows that complex conjugation induces an
isomorphism from $H_{\zeta^i}$ to  $H_{\zeta^{-i}}$.  It is also
apparent that the intersection forms on these two eigenspaces are
conjugate so that the signatures are the same.

\begin{thm} For $i \ne 0$ the map on $H_1(\tilde{F}_k,\cc)$ induced by
$\tilde{\tau}_k$ maps $H_{\zeta^i}$ to $H_{\zeta^{-i}}$ isomorphically,
and this isomorphism agrees with the one induced by complex conjugation.
\end{thm}

\begin{proof}   The action of $\tilde{\tau}_k$ and $R_k$ generate a
dihedral action on the rational vector space $H_1(\tilde{F}_k,\qq)$.  
It follows from representation theory that any such rational
representation of the dihedral group splits as the direct sum of three
types of representations
\cite{serre}. The first, denoted simply as
$\qq$, is the trivial action on $\qq$. The second, denoted $\qq_-$, is
the action on $\qq$ for which $R_k$ acts trivially and $\tilde{\tau}_k$
acts by multiplication by $-1$. The last is denoted $\qq[z]$, where $z$
is a primitive
$p$th root of unity, $R_k$ acts by multiplication by $z$, and
$\tilde{\tau}_k$ acts by conjugation, ie, $\tilde{\tau}_k(z^i) =
z^{p-i}$. 

Complexifying (ie, tensoring with $\cc$ over $\qq$) these
representations we see that the only one for which $R_k$ has a
$\zeta^i$ eigenspace is for $\qq[z]$, where that eigenspace is
generated by $1 +
\zeta^{-i}z + \zeta^{-2i}z^2 + \ldots +\zeta^i z^{p-1}$. Now the action
of $\tilde{\tau}_k$ is to conjugate; again, this map sends $z^i$ to
$z^{-i}$.  When we complexify the representation, the action of complex
conjugation fixes the $z^i$ and maps $\zeta^i$ to $\zeta^{-i}$.  
In either case we see the actions are the same. This completes the
proof.
\end{proof}

\subsection{Tangle addition of Casson--Gordon invariants} 
  
We are now in the following setting.  The knots $K$ and $K^*$ are
positive mutants, split into tangles $T^a$ and $T^b$ by a splitting
sphere $S$.  The 2--fold branched cover of $K$, $M_1$, splits as
$M_1^a \cup_{id} M_1^b$ and the 2--fold branched cover of $K^*$ splits
as
$M_1^a \cup_{\tau_1} M_1^b$ , where $\tau_1$ is an involution of the
2--fold branched cover of $S$, $F_1$.  We have a $\zz_p$--valued
character $\chi$ on $H_1(M_1)$ and a character $ \chi^*$ on
$H_1(M_1^*)$ obtained from $\chi$ by inverting it on $M_1^b$.  The goal
of this section is the proof of the following. 

\begin{thm}\label{tangleadd} In this setting,
$\sigma(K,\chi,k,i) = \sigma(K^*,\chi^*,k,i)$.\end{thm}

\begin{proof} As expected, the proof begins by decomposing the
manifolds used to compute the Casson--Gordon invariants. Hence: \vskip1ex
\noindent{\bf Step 1}\qua {\sl \ Decomposing $M_k$} 

We begin by  describing how to construct a knot $K^a$ as the union of
$T^a$ and a trivial twisted tangle, $T^c$, (see Figure 1, illustrating
a trivial twisted tangle) so that the restricted character
$\chi$ on the 2--fold branched cover of $T^a$ extends to the branched
cover of
$K^a$. There is a representation of
$\pi_1(S^3-K)$ to the dihedral group, $D_{2p} = \langle t,r | t^2 = 1,
r^p = 1, trt = r^{-1} \rangle$, so that its restriction to the index 2
subgroup gives the character $\chi$. (This observation concerning
dihedral representations is well known; we present the details in the
more complicated setting of 3--fold covers in the proof of Theorem
\ref{cohomology}.) The restriction of this representation to
$S -
\{4 \mbox{ points}\}$ extends naturally to the complement of some
standard twisted tangle, which we denote $T^c$.

\vskip.2in
\epsfxsize=1.8in
\centerline{\epsfbox{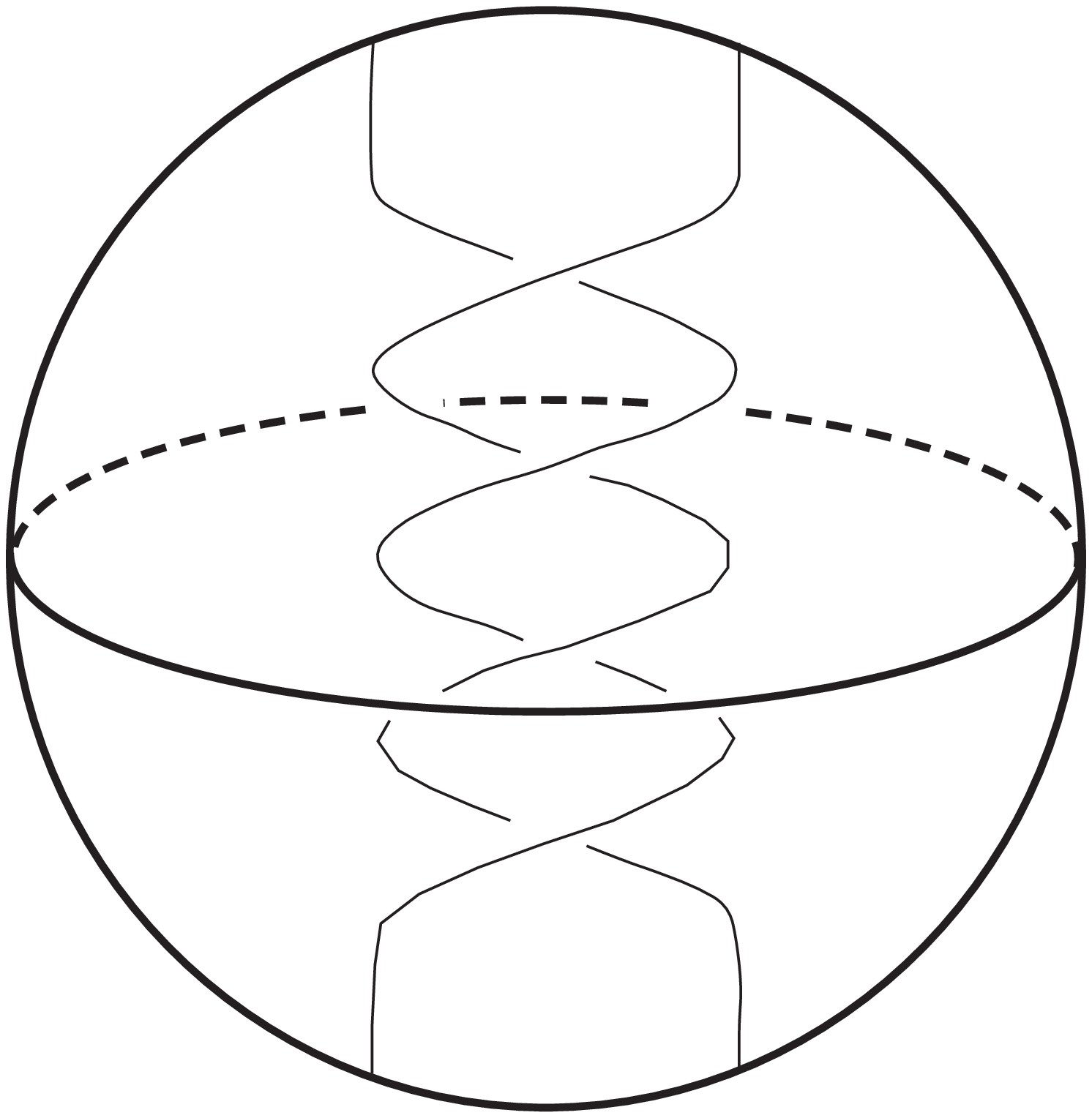}}
\vskip.1in
\centerline{\small Figure 1: A trivial twisted tangle}
\vskip.2in

Notice that the same tangle $T^c$ can be glued to the tangle
$T^b$ to yield the knot $K^b$ for which the character also extends in
the same way over the cover of $T^c$.

We let $M^a_k$, $M^b_k$, and $M^c_k$ be the induced branched covers of
$T^a$,
$T^b$, and $T^c$. Notice that $\tau$ extends to an involution of
$T^c$, and lifts to give involutions $M^c_k$.

We can form $N^a_k = M^a_k \cup_{id} M^c_k$, the $2^k$--fold branched
cover of $K^a$. Similarly, form
$N^b_k = M^b_k
\cup_{id} M^c_k$, the $2^k$--fold branched cover of $K^b$. Let
$\chi_k^a \in H^1(N_k^a;\zz_p)$ and $\chi_k^b \in H^1(N_k^b;\zz_p)$ be
the induced characters. 

\vskip2ex
\noindent{\bf Step 2}\qua {\sl \ Non-additivity of the signature} 

We now want to relate the Casson--Gordon invariants of both
$N^a_k$ and
$N^b_k$ to those of $M_k$ and $M_k^*$. As usual, we let $p$ copies of
$N^a_k$ and
$N^b_k$ bound 4--manifolds $W^a_k$ and
$W^b_k$ over $ \chi_k^a$ and $ \chi_k^a$, respectively. The desired
manifolds
$W_k$ and
$W^*_k$ needed to compute the Casson--Gordon invariants for $M_k$ and
$M_k^*$ are obtained by gluing $W^a_k$ and $W^b_k$ along
$M^c_k$ using either the identity map or the extension of
$\tau_k$ to $M^c_k$. We then want to apply additivity of the signature,
both on the level of the base manifold and for the eigenspace signature
on the
$p$--fold cyclic covers. However signatures do not add in this setting;
the necessary result is Wall's ``non-additivity'' formula for signature
\cite{wa}, which we now summarize.

Suppose that a 4--manifold $Z$ is formed as the union of $Z_1$ and
$Z_2$ along a 3--dimensional submanifold $Y$ contained in both
boundaries and the boundary of
$Y$ is the surface $X$. Then $$\mbox{sign}(Z) = \mbox{sign}(Z_1) +
\mbox{sign}(Z_2) - \psi(H_1(X),V_1,V_2,V_3),$$ where $\psi$ is a
function defined by Wall and $V_1, V_2$, and $V_3$ are the kernels in
$H_1(X)$ of the inclusion of $X$ into
$\partial Z_1 - Y$, $\partial Z_2 - Y$, and $ Y$ respectively. 

\vskip1ex
\noindent{\bf Step 3}\qua {\sl \ Applying non-additivity} 

From this most steps of the theorem follow readily. Using the notation
of  Section 2, we have $$\mbox{sign}(W_k) =
\mbox{sign}(W^a_k) + \mbox{sign}(W^b_k) -
\psi(H_1(F),V_1,V_2,V_3)$$
$$\mbox{sign}(W^*_k) = \mbox{sign}({W_k^a}^*) +
\mbox{sign}({W_k^b}^*) - \psi(H_1(F^*),V_1^*,V_2^*,V_3^*).$$
Now, by identifying $F$ (and $F^*$) with the boundary of $M^a_k$, we
see that all these terms are in fact identical, except perhaps the
terms $V_2$ and $V_2^*$. The first is the kernel of $H_1(F)
\rightarrow H_1(M^b_k)$ and the second is the kernel of this same
inclusion, precomposed with $\tau_k$. But since $\tau_k$ acts by
multiplication by $-1$, these are the same as well.

 On the eigenspace level the situation is more difficult.  First, one
needs to observe that Wall's theorem applies on the level of
eigenspaces.  At the conclusion of \cite{wa}, Wall  notes
that the results of his paper hold in the setting of $G$--manifolds and
for the $G$--signature.  Specifically, 
$$\sigma_g(W_k) = \sigma_g(W_k^a) + \sigma_g(W_k^b) -
\psi_g(H_1(\tilde{F}), V_1, V_2, V_3),$$
where $\sigma_g$ denotes the $G$--signature and $\psi_g$ is an invariant
depending only on the subspaces $V_i$ of $H_1(\tilde{F})$.  Since the
eigenspace signatures and the $G$-signatures are equivalent via a Fourier
transform by \cite{cg1,aps}, it follows that (again using
the notation of Section 2, so, for example,
$H^a_{k,i}$ means the $\zeta^i$ eigenspace  of
$H_2(\tilde{W}^a_k, \cc)$)  we have
$$\mbox{sign}(H_{k,i}) = \mbox{sign}(H^a_{k,i}) +
\mbox{sign}(H^b_{k,i}) -
\psi(H_1(\tilde{F})_{\zeta^i},V_1,V_2,V_3)$$
$$\mbox{sign}({H_{k,i}}^*) = \mbox{sign}({H^a_{k,i}}^*) +
\mbox{sign}({H^b_{k,i}}^*) -
\psi({H_1(\tilde{F})_{\zeta^i}}^*,V_1^*,V_2^*,V_3^*),$$ 
where $\psi$ is a function depending only on the subspaces $V_i$ of
$H_1(\tilde{F})_{\zeta^i}$.

Alternatively, one can check line-by-line that Wall's arguments go
through in the setting of eigenspaces and that hence 
these formulas hold.

By identifying $F_k$ with the boundary of $M^a_k$, we have that all the
terms in the sum are identical, except for two pairs that we now deal
with. 

\vskip3ex
\noindent {\bf Step 4} \ $\mbox{sign}(H^b_{k,i}) =
\mbox{sign}({H^b_{k,i}}^*)$

	These might be different because the second term is obtained after
replacing
$\chi_k$ with $-\chi_k$ on $M^b_k$. This has the effect of inverting
the $\zz_p$ action on
$\tilde{W}^b_k$, and thus interchanges the eigenspaces. However,
complex conjugation also induces an interchange of eigenspaces, and as
we noted in the previous section, this preserves signatures, so the
terms are in fact equal.
\vskip2ex
\noindent {\bf Step 5}\qua $ V_2 = V_2^*$

This is the most difficult and delicate step. We must prove that the
action of the lift of $\tau_k$, $\tilde{\tau}_k$, to the
$p$--fold cyclic cover preserves the kernels of the inclusions on the
eigenspaces. More precisely, what must be shown is that the kernel of
$\tilde{i}\!: H_1(\tilde{F}_k,\cc)_{\zeta^i}
\rightarrow H_1(\tilde{M}^b_k,\cc)_{\zeta^i}$ and the kernel of
$\tilde{i}
\circ \tilde{\tau}_k\!: H_1(\tilde{F}_k,\cc)_{\zeta^i} \rightarrow
H_1(\tilde{M}^b_k,\cc)_{\zeta^i}^*$ are the same. Here the subscript
$\zeta^i$ denotes the $\zeta^i$ eigenspace. 

As mentioned earlier, there is a natural isomorphism  from
$H_1(\tilde{M}^b_k,\cc)_{\zeta^i}$ to
$H_1(\tilde{M}^b_k,\cc)_{\zeta^{-i}}$ induced by complex conjugation and this
second eigenspace is isomorphic to
$H_1(\tilde{M}^b_k,\cc)_{\zeta^i}^*$ since the $^*$ denotes the groups
associated with
$-\chi_k$ for which the action of the deck transformation is inverted.
Also, notice that complex conjugation preserves the kernel of
inclusions, since the complex kernel is simply the real kernel tensored
with $\cc$.

In the  last section we proved that
$\tilde{\tau}_k$ induces an isomorphism of
$H_1(\tilde{F}_k,\cc)_{\zeta^i}$ to
$H_1(\tilde{F}_k,\cc)_{\zeta^{-i}}$ that agrees with the one induced by
complex conjugation. Theorem \ref{tangleadd} follows.
\end{proof}

\noindent{\bf Proof of Theorem \ref{posmut}}\qua Since  $K \# -K$ is
slice, Theorem \ref{abc} says that $K \# -K$ is cg--slice.  But, $K \#
-J$ is a positive mutant of $K \# -K$.  Hence the following theorem
completes the proof.

\begin{thm} If $K$ is cg--slice and $K^*$ is a positive mutant of $K$,
then
$K^*$ is also cg--slice.\end{thm}

\begin{proof} We have already seen that there is a correspondence between
characters on $H_1(M)$ and $H_1(M^*)$ and that, by Theorem
\ref{sigmai}, $\sigma_i$ is invariant under this correspondence. It
remains to check that metabolizers are preserved under the
correspondence. To see this we need to be look a little more closely at
the homology and the linking form. 

Via a Mayer-Vietoris argument, $$H_1(M) = (H_1(M^a) \oplus H_1(M^b)) /
G_1$$ and $$H_1(M^*) = (H_1(M^a) \oplus H_1(M^b)) / G_2,$$ where
$G_1$ is the subgroup generated by elements of the form
$(i_*(x),i_*(x)), x \in H_1(F )$, and $G_2$ is the subgroup generated
by the set of elements of the form $(i_*(x),-i_*(x)), x \in H_1(F)$. An
isomorphism is given by the map $(x,y) \rightarrow (x,-y)$.

To check that the isomorphism is an isometry of linking forms, we need
only to check that the linking number of classes of the form
$(x,0)$ or
$(0,y)$ is preserved. 

Recall that the linking form $lk(a,b)$ is defined geometrically by
letting a multiple of $a$ bound a chain and computing the intersection
number $b$ with that chain, and then dividing by the order of the
multiple.

To check that the linking number of a class of the form $(x,0)$ and a
class of the form $(0,y)$ is preserved, let $x$ bound a chain $z$ in
$M$. We can assume that $z$ is transverse to $F$ and intersects $F$ in
a $\tau$ invariant 1--chain. It follows that $z$ can then be cut and
reglued to give a chain $z*$ in
$M^*$. Notice though that since $\tau$ acts by
$-1$ on $H_1(F)$, the portion of $z*$ in
$M^b$ has its orientation reversed in the construction.  Hence, the
intersection number of $z$ with $(0,y)$ in $M$ is the same as the
intersection number of $z*$ with $(0,-y)$ in $M*$.

A similar argument applies to the cases of linking numbers of classes of
the form
$(x,0)$ and
$(x',0)$ or of classes of the form $(0,y)$ and
$(0,y')$.  

Since the isomorphism of $H_1(M)$ with $H_1(M^*)$ preserves linking, it
takes metabolizers to metabolizers.
 \end{proof}

\subsection{Examples}\label{examps1}

Let $K_a$ be the $a(a+1)$ twisted double of the unknot, illustrated in
Figure 2 below.  

\vskip.2in
\epsfxsize=2.5in
\centerline{\epsfbox{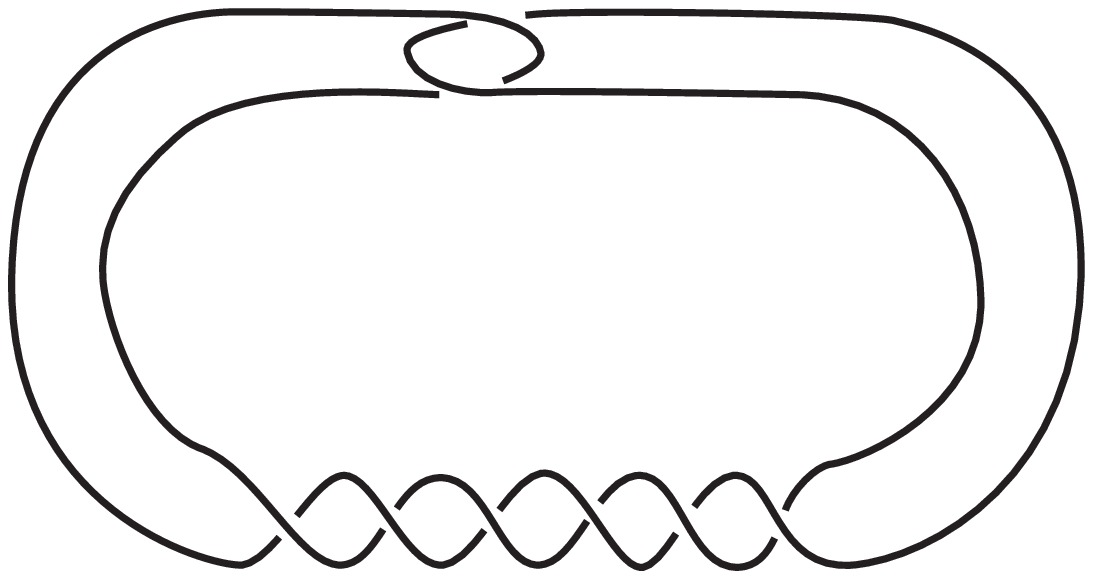}}
\vskip.1in
\centerline{\small Figure 2: A twisted double of the unknot}
\vskip.2in

 These knots were principal examples in both
\cite{cg1, cg2} and \cite{ji}. In particular, it is shown in
\cite{cg1,cg2} that for
$a>1$ these knots are algebraically slice but not slice. In \cite{ji} it
is further shown that if a sequence of values of $a>1$ are chosen so that
the corresponding sequence of values $2a+1$ are distinct primes then
the associated sequence of knots forms a linearly independent set in the
kernel of Levine's homomorphism. 

The proofs of the results of the preceding paragraph depend (in the
language developed here) on an estimate for $\sigma_b$. A careful
examination shows that the results in \cite{cg1,cg2,ji}  imply, via
Theorem \ref{posmut}, that the sequence is independent in the
kernel of
$\phi_2$. In fact one has the following theorem.

\begin{thm}\label{phi2inf}  Let $\{a_i\}_{i=1}^\infty$ be an increasing
sequence of positive integers so that the sequence 
$\{2a_i+1\}_{i=1}^\infty$ is a sequence of primes. Then the knots
$K_{a_i}$ are linearly independent in
$\calm$ and generate a $\zz^\infty$ subgroup of   $\ker\phi_2\!:\calm\to
\calg$.
\end{thm}

\begin{proof} The Seifert form for $K_a$ is given by
$$\begin{pmatrix} -1 & 1 \\ 0 & a(a+1)
\end{pmatrix}.$$ 
A symplectic basis for the form is given by the pair of
elements
$\{(a,-1) , (1,0)\}$. In this basis the Seifert form has a $0$ in the
upper left-hand entry, ie, the vector $(a,-1)$
is a metabolizing vector for the Seifert form. In particular $K_a$ is
algebraically slice, and hence lies in the kernel of $\phi$; its image
in $\calm$ lies in the kernel of
$\phi_2$.  By Corollary \ref{posmutcor} $K_a$ is non-trivial in $\calm$
if $K_a$ is not cg-slice.

 Representing the  elements $(a,-1)$ and $(1,0)$ by curves that meet
once on the Seifert surface, and taking a neighborhood of the pair, one
can describe the given Seifert surface as a disk with two bands added,
one with core the curve representing
$(a,-1)$ and the other with core representing the element
$(1,0)$. One can easily check that the first band is tied in a
$(-a,a+1)$--torus knot, which we will denote for now by $T_{-a}$.
Finally note that tying a copy of $T_a$ in the band yields a slice
knot, since it now bounds a genus one Seifert surface for which one band
has 0 framing and has a core that is slice.  We let $K^s_a$ denote
this slice knot.

The 2--fold branched cover of $K_a$ (and $K_a^s$) has homology
$\zz_{(2a+1)^2}$. If
$(2a + 1)$ is prime, there is a unique metabolizer for this form, and
hence a unique (up to multiple) character vanishing on that
metabolizer. By Theorem \ref{lith}, tying the knot in the band
changes the value of
$\sigma_b$ by the sequence
$(2^k\sigma_{j/(2a+1)}(T_a))_{k=1\ldots\infty}$. In other words, for a
character $\chi$ in the metabolizer,
$$\sigma_b(K_a,\chi, \ell)-\sigma_b(K_a^s,\chi,\ell)=
(2^k\sigma_{j/(2a+1)}(T_a))_{k=1\ldots\infty}$$
for some $j$.
The number  $j$
is determined by the value of $\chi$ on the curve refered to as
$U$ in the statement of Theorem \ref{lith}.  

A careful analysis of the
branched cover of $K_a^s$ shows that $j$ will be nonzero modulo
$2a+1$.  Details of such a calculation are presented in
\cite{gi} or \cite {gl1}.  They are based on the presentation of the
2--fold cover of $K_a^s$ using the methods developed in
\cite{ak}.

 This result must be 0  in the group of sequence modulo
bounded sequences since $K^s_a$ is slice and since the metabolizer is
unique. Therefore  the  knot $K_a$ will have
$\sigma_b(K_a,\chi,\ell) =
(2^k\sigma_{j/(2a+1)}(T_a))_{k=1\ldots\infty}$.

 A proof of the following lemma is given in the next
section.

\begin{lemma}\label{herald} Let $K$ denote the $(-a,a+1)$--torus knot
with
$a >1$. Then the signature function
$\si_t(K)$ is positive for ${1\over{a(a+1)}}<t<1- {1\over
{a(a+1)}}.$
\end{lemma}

Since $j$ is non-zero modulo $2a+1$, we may assume 
$0<j<{2a+1},$
so that  $${1\over{a(a+1)}}< \frac{j}{2a+1}<1- {1\over
{a(a+1)}}.$$ Therefore Lemma \ref{herald} implies that
$\sigma_{j/(2a+1)}(T_a)$ is  positive and so in particular $K_a$ is
not cg-slice. Thus Corollary \ref{posmutcor} implies that $K_a$ is
non-trivial in $\calm$, and hence $K_a$ is a non-trivial knot in the
kernel of $\phi_2$. 

 To construct the infinite family of linearly independent
algebraically slice knots in
$\calm$, let $\{ a_i \}_{i= 1}^\infty$ be an increasing sequence of
positive integers so that the sequence $\{ 2a_i +1\}_{i= 1}^\infty$ is
a sequence of distinct primes.  We will show the sequence
$\{K_{a_i}\}$ is linearly independent. Given a sequence of integers
$\{n_i\}$ we must show that the knot $L = \#_i n_i K_{a_i}$ is not
slice.  By reindexing, and replacing $L$ with $-L$ if necessary, we may
assume that
$n_1$ is positive.  We will show that for any character
$\chi$ on the 2--fold branched cover of $L$ taking values in
$\zz_{2a_1+1}$ the value of
$\sigma_b(L,\chi,\ell)$ is nontrivial for any $\ell$.  

Since the 2--fold branched cover of $L$ is a connected sum, we can write
$\chi$ as a sum $\#_i \chi_i$ where $\chi_i$ is the restriction of
$\chi$ to the cover of $n_iK_i$. But since  the 2--fold branched cover
of $K_i$ has no $2a_1 +1$ homology for $i>1$, we have that
$\chi_i$ is trivial for $i>1$ and therefore $\sigma_b(n_i
K_{a_i},\chi_i,\ell) = 0$ for $i >1$.
Hence by the   additivity of
$\sigma_b$, 
$$\sigma_b(L,\chi, \ell)=\sigma_b(n_1
K_{a_1},\chi_1,\ell).$$
We will show  that
$\sigma_b(n_1 K_{a_1},\chi_1,\ell)$ is nonzero.  Using additivity again
this  is the sum of values of $\sigma_b(K_{a_1},\eta_m,\ell)$ where
$\eta_m$ denotes the restriction of $\chi_1$ to the $m$th factor of
$n_1 K_{a_1} = K_{a_1}
\# K_{a_1}
\#
\cdots \# K_{a_1}$.  Explicitly
$$\sigma_b(n_1
K_{a_1},\chi_1,\ell)=\sum_{m=1}^{n_1}\sigma_b(K_{a_1},\eta_m,\ell)
=\sum_{i=m}^{n_1}(2^k\sigma_{{j_m}/(2a_1+1)}
(T_{a_1}))_{k=1,\ldots,\infty}.$$
Since $\chi$ is nontrivial, at least one of the
$\eta_m$ is nontrivial. For each such $m$,    $j_m$ is
non-zero modulo
$2a_1+1$ as before, and  Lemma \ref{herald}  implies that the
corresponding
$\sigma_{{j_m}/(2a_1+1)}(T_{a_1})$ is positive. 

For those $m$ such that $\eta_m$ is trivial,  
$\sigma_{{j_m}/(2a_1+1)}(T_{a_1})$ is zero.  Therefore,
$\sigma_b(n_1 K_{a_1},\chi_1,\ell)$ is non-zero. Since $\sigma_b(n_1
K_{a_1},\chi_1,\ell)=\sigma_b(L,\chi,\ell)$, $L$ is not cg-slice and so
represents a non-trivial element in $\calm$ by Corollary
\ref{posmutcor}. This concludes the proof of Theorem \ref{phi2inf}.
\end{proof}

\subsection{The Proof of  Lemma \ref{herald}} The proof of this result
concerning the signatures of
$(a,a+1)$--torus knots could presumably be constructed using a detailed
analysis of the Seifert form of these knots.  See
\cite{murasugi} for an example of such a computation, carried out in
the case of $1/2$ signatures in detail.    However, the recently
developed understanding of the signatures in terms of
$SU(2)$ representations offers an alternative which, though more
sophisticated in ways, brings a simple perspective.

\vskip3ex

\noindent{\bf Proof of Lemma \ref{herald}}\qua We begin by recalling that
$SU(2)$ can be identified with the unit quaternions  via
$$\begin{pmatrix} z & w \\ -\bar{w} & \bar{z}
\end{pmatrix} \leftrightarrow z+ w{\bf j}$$ for $z,w \in \cc$,
$|z|^2+|w|^2=1$. 

Results of Herald \cite{herald1} (see also \cite{heusener}) interpret
the signature function of a knot in terms of the space of conjugacy
classes of
$SU(2)$ representations of the fundamental group of the knot complement
and the trace of these representations restricted to the meridian,
$\mu$, of the knot.  Roughly stated, 
$\si_t$ is equal to  twice the algebraic count of conjugacy
classes of non-abelian $SU(2)$ representations $\rho$ of the knot group
satisfying tr$(\rho(\mu))= 2\cos(\pi t)$.  In general one must
perturb the representation space to assure transversality before
performing the count, but it can be shown that for $(p,q)$--torus knots
such a perturbation is unnecessary.  Furthermore, in general there is a
sign issue in performing the count, but it follows from the work of
Herald
\cite{herald2} that for a $(p,q)$--torus knot each representation
contributes  --sign($pq$).

It follows from the previous discussion that a proof of Lemma
\ref{herald} is implied by the statement that for all
$t$, ${1\over{a(a+1)}}<t<1- {1\over {a(a+1)}}$, 
there is a non-abelian $SU(2)$ representation, $\rho$, of the
$(-a,a+1)$--torus knot group with tr$(\rho(\mu))= 2\cos(\pi t)$. 
The rest of this proof is devoted to proving the existence of such a
representation.

Recall that $\pi=\pi_1(S^3-K)=\langle x,y\ |\ x^a=y^{a+1}\rangle$ and
that a simple calculation shows that the  meridian is represented by
$\mu =xy^{-1}$. 

Given a pair
$m,n$ of integers such that $0<m<a$, $0<n<a+1$ with 
$m \equiv n\pmod 2$, and $u\in[0,1]$, there is a  representation
$\rho_{u,m,n}\!:\pi\to SU(2)$ satisfying 

$$\rho_{u,m,n}(x) =e^{2\pi i m/2a}$$ and 
$$\rho_{u,m,n}(y) =\cos( 2\pi
 n/2(a+1)) + \sin ( 2\pi
 n/2(a+1))     (\cos(\pi u){\bf i}+\sin(\pi u){\bf j}).\eqno{(2)}$$
(These two unit quaternions are easily shown to be $a$ and $a+1$ roots
of
$1$ or
$-1$, depending on whether $m$ and $n$ are even or odd.)  The
representation $\rho_{u,m,n}$ is nonabelian unless $u = 0$ or 1.

Setting $u = 0$  and $1$ gives $$\rho_{0,m,n}(\mu) =
\exp(2\pi i({m\over{2a}}-{n\over{2(a+1)}}))$$ and $$\rho_{1,m,n}(\mu) =
\exp(2\pi i({m\over{2a}}+{n\over{2(a+1)}})).$$
The traces of these values, given as  twice the real parts of the unit
quaternions, are $2\cos( \pi  ({m\over{a}}-{n\over{(a+1)}}))$  and
$2\cos( \pi  ({m\over{a}}+{n\over{(a+1)}}))$.  By continuity (in $u$)
all traces between these two values occur also.  

That the corresponding open intervals cover
$({1\over{a(a+1)}}, 1-{1\over{a(a+1)}})$ can be seen just by
considering the pairs $(m,n)= (1, 2k+1), \ k=0,1, 2,\cdots$.
\endproof

\section{The kernel of $\phi_2$ contains $\zz_2^\infty$}

In \cite{li3} it was shown that the kernel of Levine's homomorphism
contains a subgroup isomorphic to $\zz_2^\infty$.  Here we will show
that those examples remain nontrivial in $\calm$ and that the kernel of
$\phi_2$ contains an infinite collection of elements of order 2.  The
proof is much like that of \cite{li3} and we will only outline the
arguments, highlighting the points where they have to be enhanced.

Let $K_T$ be the knot illustrated in Figure 3.  The figure illustrates
a genus 1 Seifert surface for $K_T$ in which one band has the knot $T$
tied in it and the other has $-T$ tied in it.  The bands are twisted so
that the Seifert form is
$$\begin{pmatrix} 1 & 1 \\ 0 & -1
\end{pmatrix}.$$ In the case that $T$ is the unknot $K_T$ is the Figure
8 knot, $K_0$.  In general it is easily seen that $K_T = -K_T$ and
hence that $K_T$ is of order 1 or 2 in  $\calc$.  (Fox and Milnor were
the first to show that the order of $K_0$ is exactly 2, using the
Alexander polynomial.) 

\vskip.2in
\epsfxsize=2.5in
\centerline{\epsfbox{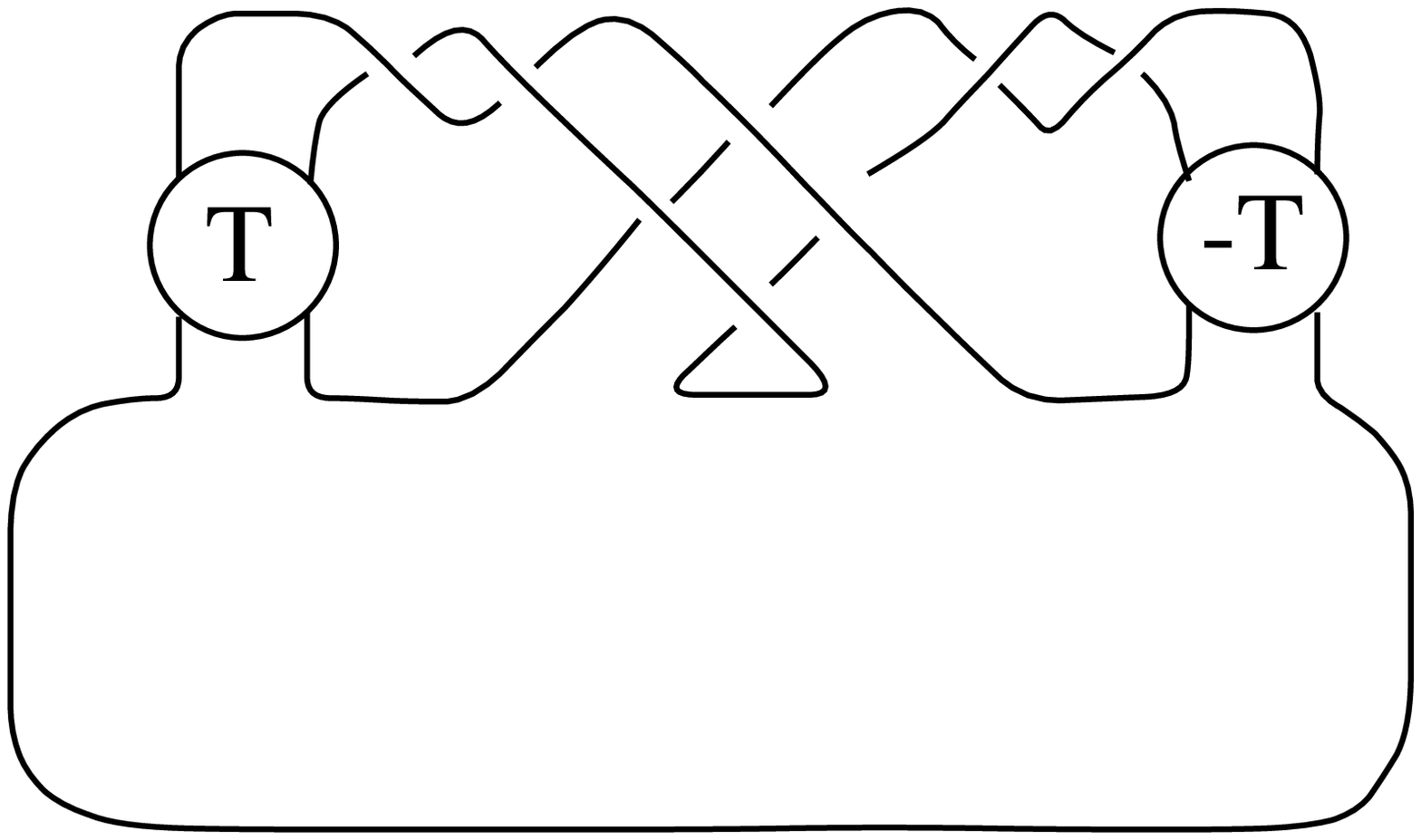}}
\vskip.1in
\centerline{\small Figure 3: An order 2 knot}
\vskip.2in

The examples we want to consider are knots of the form $ K_T \#
K_0$, for particular knots $T$.  Note immediately that these knots are
of order either 1 or 2 in $\calm$.  Also, they are in the kernel of
$\phi_2$, since the two summands have the same Seifert form, the
Seifert form of an order 2 knot.  We will  show that for an
appropriate set of knots $\{T_i\}_{i=1}^\infty$, the knots $K_{T_i} \#
K_0$
are distinct in $\calm$.  Hence,
the kernel of
$\phi_2\!:\calm\to\calg$
contains $\zz^\infty_2$. 

Showing that 
$K_{T_i}
\# K_0 \ne K_{T_j} \#
K_0$ in $\calm$ is clearly equivalent to showing that
$K_{T_i} \#
K_0 \# K_{T_j} \#
K_0 $ is non zero in $\calm$.  But this connected sum is the same (up to
concordance) as $K_{T_i} \# K_{T_j}$, since the $K_0$ summands cancel.

\begin{thm}\label{z2inf}  For the appropriate choice of
$\{T_i\}_{i=1}^\infty$, the knot $K_{T_i} \# K_{T_j}$ is not
cg--slice, for $i \ne j$.\end{thm}

\begin{proof} The argument begins with some basic computations that follow
readily from the techniques of \cite{rol} and \cite {ak}.  Details are
presented in
\cite{li3}.  The knot $K_T$ is built from $K_0$ by removing
neighborhoods of curves $B_1$ and $B_2$ linking the bands and replacing
them with the complements of $T$ and $-T$.  The 2--fold branched cover
of
$K_T$, $M(K_T)$, has first homology $\zz_5$, and if we consider a map of
$M(K_T)$ to $\zz_5$ taking value $1$ on a lift of $B_1$, it takes value
$-1$ on the other lift, and takes values $2$ and $-2$ on the lifts of
$B_2$.

The homology of the 2--fold cover of $K_{T_i} \# K_{T_j}$ is $\zz_5
\oplus
\zz_5$, and a careful examination of the linking form reveals that for
any metabolizer, one of the characters that vanishes on that metabolizer
will take value 1 on a generator of the first  $\zz_5$ summand and $\pm
2$ on the corresponding generator of the other summand.  (In brief, a
metabolizing element
$(a,b)$ must satisfy $a^2 + b^2 = 0$ mod 5.)

Hence, a calculation similar to the one for doubled knots given in the
previous section shows that for this character
$$\sigma_b(K_{T_i} \# K_{T_j}) = \left(  2^k 
\left( \sigma_{1/5}(T_i) -  \sigma_{2/5}(T_i)     +\sigma_{2/5}(T_j) - 
\sigma_{1/5}(T_k)            \right)
\right)_{k=1\ldots
\infty}$$
An explicit calculation for the $(2,7)$--torus knot, $T_1$, gives that
$\sigma_{1/5}(T_1) = -2$ and $\sigma_{2/5}(T_1) = -6$.  Hence, if we let
$T_i$ denote the connected sum of  $i$ $(2,7)$--torus knots, one
computes 
$$\sigma_b(K_{T_i} \# K_{T_j}) = \left(  2^k 
\left( -2i + 6i     -6j  +2j            \right)
\right)_{k=1\ldots
\infty}  =  
\left(  2^{k+2} 
\left( i    -j            \right)
\right)_{k=1\ldots
\infty}$$ Clearly, this last sequence is unbounded unless $i = j$.
\end{proof} 

 Theorem \ref{z2inf} and
Corollary  \ref{posmutcor}  immediately imply the following.

\begin{cor} \label{2torinf} For the appropriate choice of knots
$\{T_i\}_{i=1}^\infty$, the knots $K_{T_i}\#K_0$ generate an infinite
2-torsion subgroup of $\ker
\phi_2\!:\calm\to \calg$. 
\end{cor}

\section{The kernel of $\phi_1$  contains $\zz^\infty$}
 
In Section 5.1 we will describe the construction of a knot, $K_J$, based
on an arbitrary knot $J$. It follows from the construction that
$K_J$ is slice; our focus will be on the positive mutants,
$K_J^*$.  For instance, Theorem 5.2 states that for some $J$, $K_J^*$
is not slice; this proves that the kernel of  $\phi_1$ is nontrivial.

In order to apply Theorem 2.3 to $K_J^*$, and linear combinations of
$K_J$, we need to understand the homology of the appropriate branched
covers.  The 3--fold  cover is sufficient for our needs, and in Section
5.2 we discuss the homology of this space. This analysis can be
done by considering the cover directly, but we will instead exploit a
relationship between the homology of the 3--fold branched cover and
metacyclic representations of the knot group.  The advantage of this
approach is that it simultaneously gives us information about the
iterated covers used in defining the Casson--Gordon invariants.

Section 5.3 presents the proof that for appropriate $J$, $K_J^*$ is not
slice.  In this section we also develop many of the basic tools and
computations needed for the linear independence results of the
following sections.

The last three sections, Sections 5.4, 5.5, and 5.6, are devoted to
proving that the kernel contains $\zz^\infty$, generated by knots
$K_J^*$ for appropriate $J$.  Most of the work is focused on
understanding the invariant metabolizers in the homology of the 3--fold
branched covering spaces.  If we consider the connected sum of $n$ of
the
$K_J^*$, the homology of the 3--fold branched cover, working with
$\zz_7$ coefficients, is $(\zz_7
\oplus \zz_7)^n$.  Under the $\zz_3$ action, this will be seen to split
as a 2--eigenspace and a 4--eigenspace, each isomorphic to 
$(\zz_7)^n$.  Metabolizers will similarly split, say for now as $M_2$
and $M_4$.  The calculations of Section 5.3 will quickly imply that if
either $M_2$ and $M_4$ contains an {\sl odd} vector, meaning one with
an odd number of nonzero coordinates, then the corresponding knot will
not be slice.  Simple linear algebra will show that if the metabolizer
is not ``evenly split'' between $M_2$ and $M_4$ then this will be the
case.  The final case, in which the metabolizer contains no odd vectors
is the most delicate.  It calls on a careful examination of the linking
form of the 3--fold branched cover, carried out in Section 5.5, along
with some difficult linear algebra, carried out in Section 5.6.

\subsection{Building the examples, $K_J^*$}

  Figure 4 illustrates the needed components of the construction of our
examples. The figure includes a knot in
$S^3$,
$K$, and three unknotted circles, $B_1$,
$B_2$ and
$B_2^*$, in the complement of $K$. Also drawn are dotted lines,
$\alpha_1$ and
$\alpha_2$. Not drawn is a second knot, $J$, to be used throughout the
construction. The final knots we construct will depend on the choice of
$J$ and will be denoted $K_J$.

\vskip.2in
\epsfxsize=3.2in
\centerline{\epsfbox{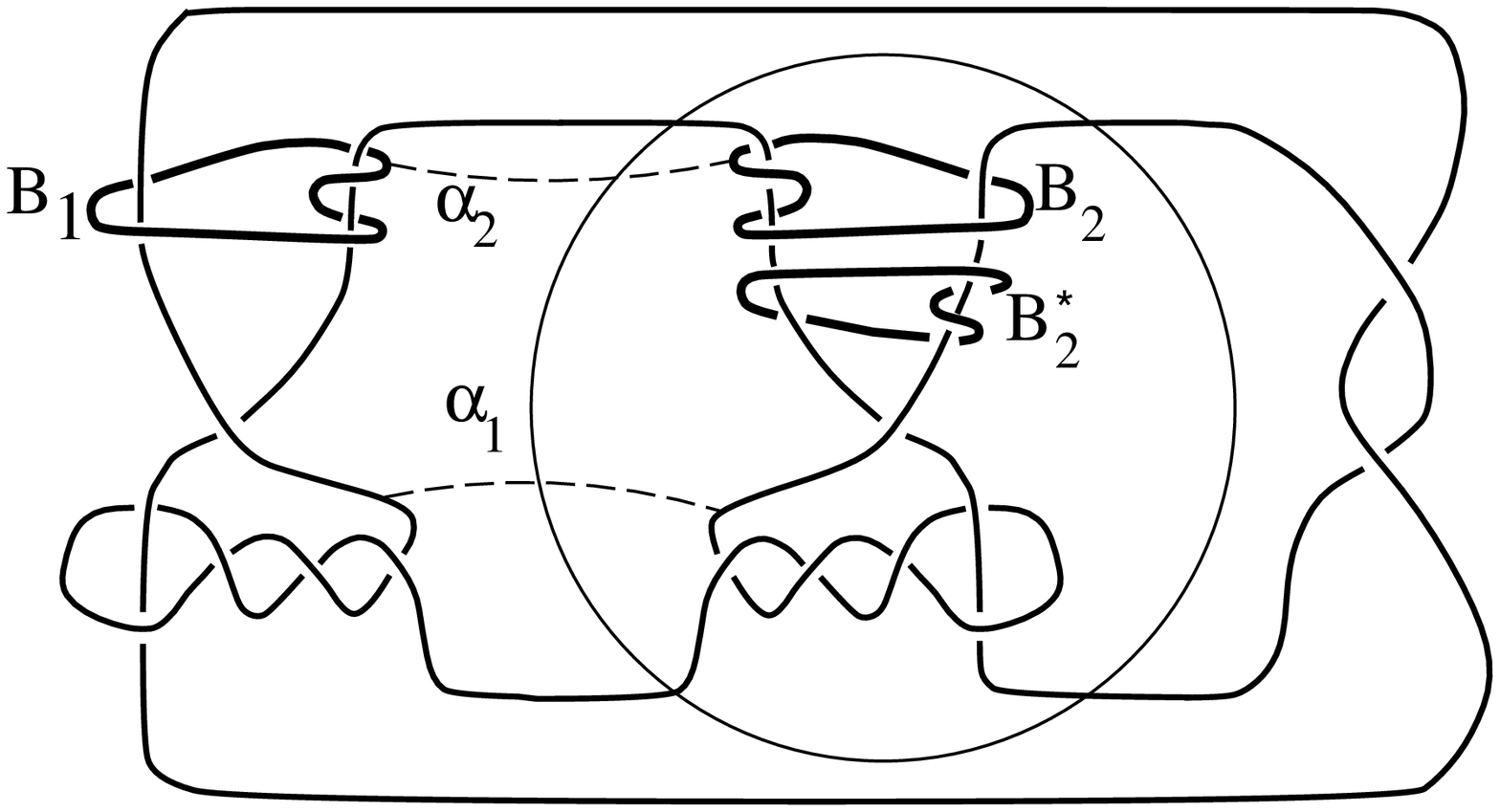}}
\vskip.1in
\centerline{\small Figure 4: The Knot $K$}

\vskip.2in

First we note that $K$ is a slice knot. This is seen by performing the
band move to
$K$ along the arc $\alpha_1$. A pair of unknotted and unlinked circles
results, and so
$K$ is slice.

The knot $K_J$ is formed by removing neighborhoods of $B_1$ and $B_2$
and replacing them with the complement of $J$ and $-J$, respectively. 
As usual, the meridian and longitude are interchanged.  As explained in
Section 2, the knot  $K_J$ is an (iterated) satellite knot in
$S^3$ since $B_1$ and $B_2$ are unknotted. We have:
 
\begin{thm} $K_J$ is slice.\end{thm}
\begin{proof} First we note that $B_1$ and $B_2$ cobound an annulus $A$
embedded in the complement in $B^4$ of the slice disk for $K$ just
constructed. To see this,  perform the band move joining $B_1$ and
$B_2$ along the arc
$\alpha_2$. The resulting circle is unknotted and unlinked from the
pair of circles formed during the slicing of $K$, so that circle can be
capped off with a disk in the complement of the slice disk to yield
$A$. We also need to note that $A$ is trivial when viewed as an annulus
in $B^4$. That is, the pair $(B^4,A)$ is equivalent to the pair
$(B^3,U) \times [0,1]$, where $U$ is the unknot. This is easily seen by
ignoring $K$ and the slice disk during the construction of $A$. 

It is clear that if one removes from $(B^3,U) \times [0,1]$ a
neighborhood of $U
\times [0,1]$ and replaces it with $(S^3 - J) \times [0,1]$ then
$B^3 \times [0,1] = B^4$ results. The boundary, $S^3 = \partial B^4$,
is obtained from the original
$S^3$ by removing neighborhoods of $B_1$ and
$B_2$ and replacing them with the complement of $J$ and $-J$
respectively, as desired. During this construction, the slice disk for
$K$ becomes the slice disk for $K_J$.\end{proof}

\vskip2ex
\noindent{\bf The positive mutant of $K_J$} 

A positive mutant of $K$ (or $K_J$) is formed by cutting along the
sphere $S$, illustrated by the circle in Figure 4, and rotating 180
degrees about the vertical axis. It is not hard to show that the
resulting mutant,
$K^*$, of
$K$  is isotopic to $K$. However, the mutant $K^*_J$ of
$K_J$ is (as we will show) not concordant to $K_J$ for many choices of
$J$. For now it should be apparent that
$K_J^*$ is formed from $K$ by removing the unknots $B_1$ and $B_2^*$
and replacing them with the complements of $J$ and $-J$ respectively.
We will ultimately prove:

\begin{thm} \label{KJnotslice}For appropriate choices of $J$, $K_J^*$
is not slice, and in particular
$K_J^*$ is not concordant to $K_J$. \end{thm} 

By considering various choices for $J$ we will construct the infinite
family of linear independent examples desired. The main result of this
section is the following.

\begin{thm}\label{phi1ker}There exists an infinite collection of knots
$J_1, J_2,\ldots$ so that for any choice of integers $n_1,n_2,\ldots$
with only  finitely many of the $n_i$  nonzero, the connected sum
$$\#_i\ n_iK_{J_i}$$ is slice, but
$$\#_i\ n_iK_{J_i}^*$$ is not slice. In particular the kernel of
$\phi_1\!:{\cal C}\to {\cal M}$ contains a subgroup isomorphic to $
\zz^\infty$.
\end{thm}

\subsection{The homology of the 3--fold cyclic cover of $K_J^*$} The
details concerning the homology of the 3--fold cyclic branched cover of
$K$ that we need are captured by the following theorem. Isolating out
this detailed information permits us to avoid a complete analysis of
the homology and cohomology, with group action, of the 3--fold cover.

\begin{thm} \label{cohomology}

\begin{enumerate}
\item[]
\item [\rm(a)] The 3--fold cover of $S^3$ branched
over $K$, $M_3$, satisfies
$H_1(M_3;\zz)=\zz_{49}\oplus \zz_{49}$. 

\item [\rm(b)] The group $H^1(M_3;\zz_{49})\cong\zz_{49}\oplus
\zz_{49}$ splits as the direct sum of two 1--dimensional eigenspaces,
$C_{18}$ and $C_{30}$, under the action of the deck transformation
$T$. (The eigenvalues are $18$ and $30$, the two nontrivial cube roots
of unity in $\zz_{49}$.)

\item [\rm(c)] There exists a choice of lifts, $\tilde{B}_1,\
\tilde{B}_2, \tilde{B}_2^*$, of the curves $B_1 , B_2,$ and $B_2^*$,
and characters  $\chi_{18} \in C_{18}$ and $\chi_{30} \in C_{30}$ (which
generate these eigenspaces), so that $\chi_i(\tilde{B}_1) = 1$,
$\chi_i(\tilde{B}_2) = 1$, and
$\chi_i(\tilde{B}_2^*) = -1$ for $i = 18, 30$.
\end{enumerate}
\end{thm}

\begin{proof} The Alexander polynomial of $K$ is
$(2t^2-5t+2)^2$. This can be quickly seen by noticing that changing one
of the two rightmost crossings on $K$ in Figure 4 changes $K$ into the
connected sum of two copies of the 3--twisted double of the unknot,
which has Alexander polynomial $2t^2-5t+2$. The Conway skein relation
shows that changing this crossing does not change the Alexander
polynomial, since the link obtained by smoothing the crossing is a
unlink. It follows that the order of $H_1(M_3;\zz)$ equals
$(2\zeta^2-5\zeta+2)^2(2(\zeta^2)^2-5(\zeta^2)+2)^2=7^4$ (see, for
example, \cite[Theorem 8.21]{bz}). The assertion that
$H_1(M_3;\zz)=\zz_{49}\oplus\zz_{49}$ will follow once we show that
$H^1(M_3;\zz_{49})=\zz_{49}\oplus\zz_{49}$.

The cube roots of $1$ in $\zz_{49}$ are $1,18$ and $30$. The deck
transformation
$T$ acting on
$H^1(M_3;\zz_{49})$ satisfies $T^3 - 1 = 0$. Hence,
$(T-1)(T-18)(T-30) = 0$. 
If $ H^1(M_3;\zz_{49})^{\zz_3}$ denotes the
fixed cohomology, then the composite
$$H^1(M_3;\zz_{49})^{\zz_3}\xrightarrow{\tau^*}H^1(S^3;\zz_{49})
\xrightarrow{\pi^*}H^1(M_3;\zz_{49})^{\zz_3}$$
is just multiplication by $3$, where $\tau^*$ denotes the transfer
and $\pi\!:M_3\to S^3$ is the branched cover (see eg \cite[Section
III.2]{bredon}); hence $H^1(M_3;\zz_{49})^{\zz_3} = 0$.
Thus $T-1$ is injective on the finite group
$H^1(M_3;\zz_{49})$ and hence an isomorphism.  It follows that 
$(T-18)(T-30) = 0$. 

Let $$C_{18}=\ker (T-18)\!:H^1(M_3;\zz_{49})\to H^1(M_3;\zz_{49})$$ and
$$C_{30}=\ker (T-30)\!:H^1(M_3;\zz_{49})\to H^1(M_3;\zz_{49}).$$ For any
$x\in H^1(M_3;\zz_{49})$,
$x = (T-30)(4x) + (T-18)(-4x)$ so
$H^1(M_3;\zz_{49})= C_{18}+C_{30}$. If $x\in C_{18}\cap C_{30}$ then
$30x=18x$, so $22x=0$ which implies $x=0$ since $22$ is relatively prime
to $49$.

Hence, the cohomology does split as the direct sum of an 18-- and
a 30--eigenspace,
$H^1(M_3;\zz_{49})=C_{18}\oplus C_{30}$. We will show that these are
each 1--dimensional (ie, cyclic of order $49$),
by identifying elements in each eigenspace with certain metabelian
representations of
$\pi_1(S^3-K)$.

Define the metacyclic group
$$G=
\langle t,r | t^3 = 1, r^{49} =1, tr= r^{30}t\rangle.$$ Consider the
elements of
$H^1(M_3;\zz_{49})$ as homomorphisms $\chi\!:\pi_1(M_3)\to \zz_{49}$. For
each $\chi\in C_{30}$, we construct a representation
$\bar{\chi}\!:\pi_1(S^3 - K)\to G$ taking value $ r^it$ on meridians, as
follows. 

Pick a basepoint for $\pi_1(S^3 - K)$ on the boundary of a meridinal
disk for $K$ and choose the basepoint of $M_3$ to lie above it, and
hence on the boundary of a meridinal disk for the branch set in $M_3$.
We can then lift a loop $\alpha\in \pi_1(S^3-K)$ to a loop
$\tilde{\alpha}$ in $\pi_1(M_3)$ by letting $\tilde{\alpha}$ be the
lift of $\alpha$ followed (if necessary) by an arc back to the
basepoint in the meridinal disk in $M_3$. Let $|\cdot |\!:\pi_1(S^3-K)\to
\zz_3$ be the homomorphism taking the meridian of $K$ to
$1\in\zz_3$.

To a class
$\chi
\in C_{30}$ we assign a representation
$\bar{\chi}$ on $\pi_1(S^3 - K)$ by the formula $$\bar{\chi}(\alpha) =
r^{\chi(\tilde{\alpha})}t^{|\alpha|}.$$ 
If $\mu\in\pi_1(S^3-K)$ denotes the meridian of
$K$ then
$\bar{\chi}(\mu^3)=1$. The projection induces an inclusion
$\pi_1(M_3)\to
\pi_1(S^3-K)/\langle \mu^3=1\rangle$ and so we can recover $\chi$ by
restricting
$\bar{\chi}$ to $\pi_1(M_3)$; its image lies in the $\zz_{49}$ subgroup
of $G$ generated by $r$.

To turn this into a one-to-one correspondence, one needs to identify
representations of
$\pi_1(S^3 - K)$ to $G$ that are equivalent on the $\zz_{49}$ subgroup.
These are precisely the conjugacy classes of representations.

Figure 5 below describes a representation of $\pi_1(S^3-K)$ to $G$. The
edges of the knot diagram (ie, the Wirtinger generators of
$\pi_1(S^3-K)$) are labeled with numbers  in $
\zz_{49}$, where the label ``$b$'' means the corresponding Wirtinger
generator is sent to the element $r^b t$ in $G$. To give a well defined
representation, the Wirtinger relation
$x_ix_jx_i^{-1}=x_k$ at a crossing implies that the relation
$r^{b_i}t r^{b_j}t (r^{b_i}t)^{-1}=r^{b_k}t$ must hold, where
$b_{\ell}$ denotes the label on the $\ell$th Wirtinger generator. This
relation is satisfied if and only if the labels satisfy $29b_i -30 b_j
+ b_k = 0$ at each crossing.

Observe that conjugation by $r$ takes $r^it$ to $r^{i-1}t$ and so in
enumerating conjugacy classes of representations (and hence
$C_{30}$), we can assume that any one label of the diagram is $0$. Or,
put otherwise, two labelings determine conjugate representations if
all corresponding labels differ by the same amount. It is a simple
exercise to check that all labelings are obtained from the labeling
of Figure 5 by either multiplying each label by $a\in\zz_{49}$ or by
adding the same number to each.

It follows that $C_{30}$ is isomorphic to $\zz_{49}$, generated by the
character $\chi$ determined by the labeling in Figure 5. A similar
computation applies for
$C_{18}$. Hence
$H^1(M_3;\zz_{49})=\zz_{49}\oplus
\zz_{49}$.

\begin{figure}[ht!]
\epsfxsize=3.7in
\centerline{\epsfbox{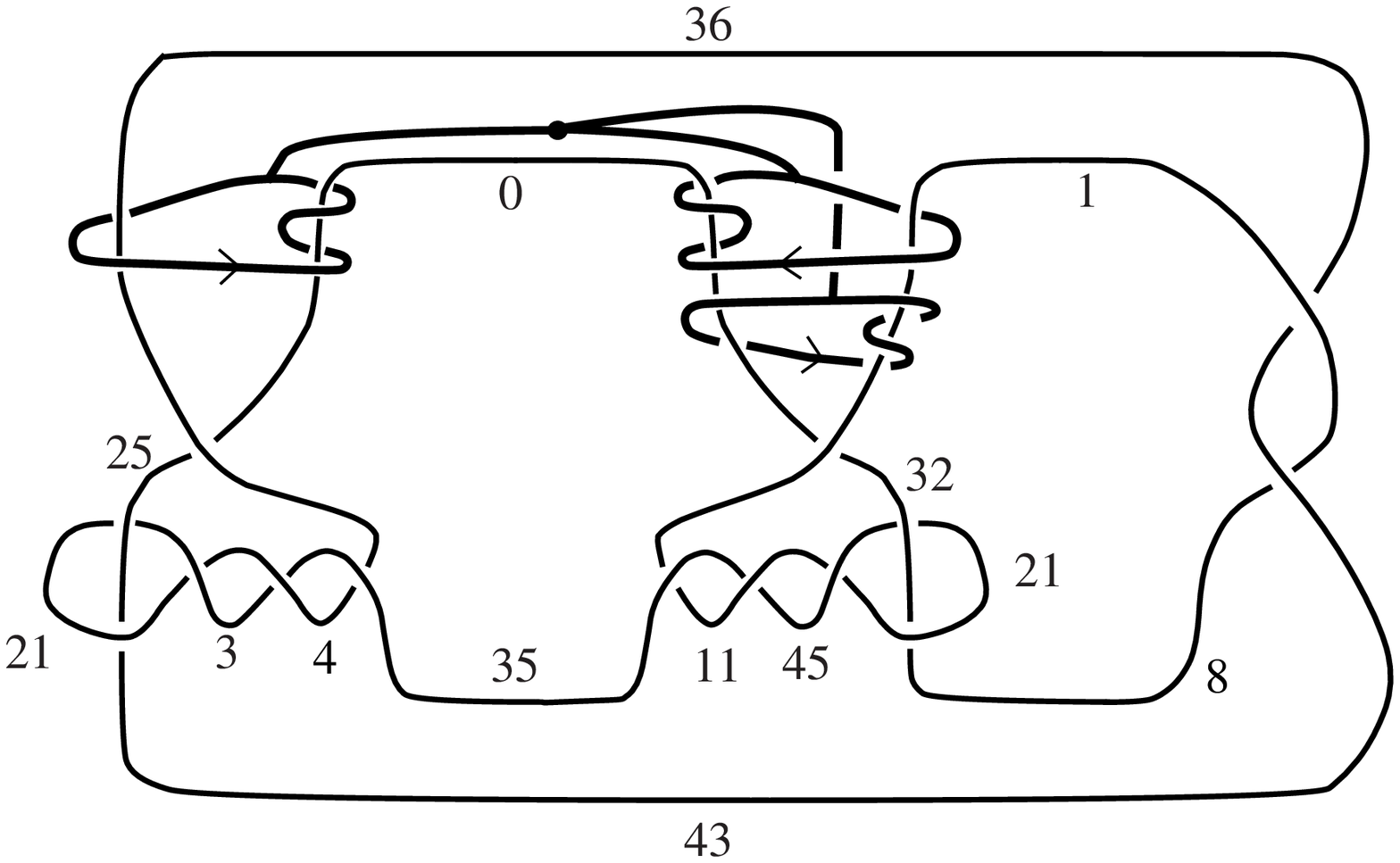}}
\vskip.1in
\centerline{\small Figure 5: The representation $\bar{\chi}\!:\pi_1(S^3-K)\to
G$}
\end{figure}

To prove part (c) one needs only compute the value of the 
representation
$\bar{\chi}$ on the
$B_i$. In Figure 5 we have chosen based loops to represent the curves
$B_1, B_2,$ and
$B_2^*$. The basing determines specific lifts $\tilde{B}_1,
\tilde{B}_2$, and
$\tilde{B}_2^*$ of these curves to $M_3$.

For the representation given by the labeling in  Figure 5 it is easy
to compute that
$\bar{\chi}(B_1)= r$, $\bar{\chi}(B_2)= r$, and $\bar{\chi}(B_2^*)=
r^{-1}$. Hence the corresponding character $\chi_{30}\in C_{30}\subset
H^1(M_3;\zz_{49})$ satisfies
$\chi_{30}(\tilde{B}_1)=1$, $\chi_{30}(\tilde{B}_2)=1$, and
$\chi_{30}(\tilde{B}_2^*)=-1$. A similar argument applies to
$\chi_{18}\in C_{18}$.
\end{proof}

Given a knot $J$ in $S^3$, we will use the notation $M_3(J)$ to denote
the 3--fold branched cover of $J$.

\begin{cor} The 3--fold branched cover $M_3(K_J^*)$ of $K_J^*$ has the
same homology and cohomology as as the 3--fold branched cover
$M_3=M_3(K)$ of
$K$; there is a natural correspondence between characters on $M_3(K)$
and characters on $M_3(K^*_J)$. In particular,  the characters
$\chi_i$ of Theorem \ref{cohomology}  extend to $M_3(K_J^*)$ and
$H^1(M_3(K_J^*);\zz_{49})=C_{18}\oplus C_{30}$. \end{cor}
\begin{proof} The curves $B_1,B_2^*$ each lift to three curves in $M_3$,
and the Mayer-Vietoris sequence shows that replacing the six solid tori
lying above these curves with  knot complements, $S^3-J$ or $S^3-
(-J)$, does not affect the cohomology.\end{proof} 

For most of our work we will only need to work modulo $7$ instead of
modulo $49$. Thus we consider the inclusion $\zz_{7}\to \zz_{49}$
induced by $1\mapsto 7$. This induces an injection
$H^1(M_3;\zz_{7})\to H^1(M_3;\zz_{49})$ with image
$7H^1(M_3;\zz_{49})$. We will always identify $H^1(M_3;\zz_{7})$ with
$7H^1(M_3;\zz_{49})$ in this way. Let
$C_2=C_{30}\cap H^1(M_3;\zz_{7})$ and let $C_4=C_{18}\cap
H^1(M_3;\zz_{7})$. Then $C_2$ and $C_4$ are $2$ and
$4$ eigenspaces of the action of $\zz_3$ on $H^1(M_3;\zz_{7})$; each is
1--dimensional, and $H^1(M_3;\zz_{7})=C_2\oplus C_4$. We let
$\chi_2=7\chi_{30}$ and $\chi_4=7\chi_{18}$; thus $\chi_2$ generates
$C_2$ and $\chi_4$ generates $C_4$ and
$\chi_i(B_1)=1=\chi_i(B_2)$ and $\chi_i(B^*_2)=-1$ (in $\zz_7$) for
$i=2,4$.  Until Section \ref{kerphi1} we will work mod 7.

\subsection{$K_J^*$ is not slice} The proof that the kernel of
$\phi_1$ is infinitely generated will depend on constructing an
infinite family of knots from $K$ by using different choices of $J$. To
introduce some of the key ingredients we will first prove Theorem
\ref{KJnotslice}, stating that for an  appropriate choice of $J$,
$K_J^*$ is not slice. 

\noindent{\bf Proof of Theorem \ref{KJnotslice}}\qua We will determine the
needed properties of
$J$ in the course of the argument.

If $K_J^*$ is slice, by Theorem \ref{cassongordon} there is some
metabolizer
$A
\subset H_1(M_3(K_J^*))$ which is invariant under $T$ and for which the
appropriate Casson--Gordon invariants (certain $\delta$) will vanish.
Simply by order considerations there will be a nontrivial
$\zz_7$ valued character, $\chi \in H^1(M_3(K_J^*);\zz_7)$, that
vanishes on
$A$. (The quotient, $H_1(M_3(K_J^*))/ A$ contains 7--torsion and hence
admits a surjection to $\zz_7$.)  The set
$A^*$ of all elements in
$H^1(M_3(K_J^*);\zz_7)$ that vanish on
$A$ is $T$--invariant since $A$ is, so $A^*$ either is 1--dimensional
and generated by an eigen-character, or $A^*$ is all of
$H^1(M_3(K_J^*);\zz_7)$. In either case, one of $\chi_2$ or $\chi_4$
must be in
$A^*$. We denote it simply by $\chi$ below.

The Casson--Gordon invariant $\delta(K_J^*,\chi)$ associated to
$\chi$ was defined in Section \ref{cgsection}. Notice that $M_3(K_J^*)$
decomposes into 7 pieces, the 3--fold cover branched cover of $S^3-(
B_1\cup B_2^*)$ branched over $K$, 3 path components homeomorphic to
$S^3-J$ lying over $B_1$, and 3 path components homeomorphic to
$S^3-(-J)$ lying over
$B_2^*$. 

We claim that
$$\delta(K_J^*,\chi)=\delta(K,\chi) \Delta_J(\zeta t)\Delta_J(\zeta^2
t)\Delta_J(\zeta^4 t)\Delta_J(\zeta^{-1} t)\Delta_J(\zeta^{-2}
t)\Delta_J(\zeta^{-4} t).$$ where

\begin{enumerate}
\item$\Delta_J$ is the Alexander polynomial of $J$.
\item$\zeta$ is a primitive 7--root of unity.
\end{enumerate}

This formula follows from several applications of Theorem
\ref{deltathm} and the following observation. If
$\chi$ takes value 1 on some lift of
$B_1$, it takes values 2 and 4 on the other lifts (if $\chi\in C_2$ then
$\chi(T\tilde{B}_1)=2\chi(\tilde{B}_1)=2$ and
$\chi(T^2\tilde{B}_1)=4\chi(\tilde{B}_1)=4$; similarly for $\chi\in
C_4$, since
$4^2=2\pmod{7}$). This explains the appearance of the terms
$\Delta_J(\zeta t),
\Delta_J(\zeta^2 t), \Delta_j(\zeta^4 t)$. We argue similarly for
$B_2^*$ using the fact that the Alexander polynomial of $-J$ is the
same as the polynomial for
$J$ and that $\chi$ takes the value $-1$ on $\tilde{B}_2^*$. 

There is a tricky point here. If the orientation of $B_1$ is changed
then $\chi$ will take value $-1$ on the lift of $B_1$. However, the
representation of the homology of the cyclic cover of $S^3 -K$ to
$\zz$ used to define $\delta$ will take value $-1 
\in \zz$ rather than $1$. Hence the term $\Delta_J(\zeta t)$ will be
replaced with
$\Delta_J(\zeta^{-1}t^{-1})$. But this just equals $\Delta_J((\zeta
t)^{-1})$, and by the symmetry of the Alexander polynomial the term is
unchanged. 

Our first conditions on $\Delta_J$ (and hence on $J$) is that it be
irreducible in
$\qq(\zeta)[t,t^{-1}]$. This is easily attained: let $J$ be a knot
having Alexander polynomial $at^2 +(1-2a)t +a$. This is irreducible
since the roots are complex, and as long as its discriminant $(1-4a)$ is
not divisible by 7 it will not factor in
$\qq(\zeta)[t,t^{-1}]$ (this is shown in greater generality in the next
section).

A second condition on $\Delta_J$ is that it be relatively prime to
$\delta(K,\chi)$ in
$\qq(\zeta)[t,t^{-1}]$ for each $\chi\in H^1(M_3;\zz_7)$. But since
$\delta(K,\chi)$ has only a finite number of irreducible factors, and
$H^1(M_3;\zz_7)$ has 49 elements, this too is easy to achieve.

With these conditions, it is clear that $\delta(K_J^*,\chi)$ can be a
norm (that is, of the form $g(t)\overline{g(t)}$) only if $$S (t) =
\Delta_J(\zeta t)\Delta_J(\zeta^2 t)\Delta_J(\zeta^4 t)\Delta_J(\zeta^{-1}
t)\Delta_J(\zeta^{-2} t)\Delta_J(\zeta^{-4} t)$$ is a norm.

We next note that the six factors in this product are irreducible and
are pairwise relatively prime in
$\qq(\zeta)[t,t^{-1}]$. Certainly the change of variable cannot make
them reducible. If two were not relatively prime they would be
associates; one would be an nontrivial $\qq(\zeta)$ multiple of the
other. But this is immediately checked to not be the case.

So, suppose that $S (t)$ is of the form $g(t)\overline{g(t)}$. Notice
that
$\overline{\Delta(\zeta^i t)} = \Delta(\zeta^{-i}t^{-1}) =
\Delta((\zeta^i t)^{-1})$. However, by the symmetry of the Alexander
polynomial, this last term is simply
$\Delta(\zeta^i t)$ (modulo a multiple of the form $t^2$ which is
itself a norm). It follows that for $S (t)$ to be a norm, it would have
to be a square; since each irreducible factor is distinct, this is not
the case. Thus $K_J^*$ is not slice.
\endproof

We formalize the last paragraph of this proof in the following useful
algebraic lemma.

\begin{lemma}\label{squares} Let $q_1(t), q_2(t),\ldots, q_n(t)$ be
integer Laurent polynomials. Let $\zeta=e^{2\pi i/7}$. Suppose that
$q_i(t)=q_i(t^{-1})$, that each $q_i(t)$ is   irreducible in
$\qq(\zeta)[t,t^{-1}]$, and that the $q_i(t)$ are relatively prime in
$\qq[t,t^{-1}]$.  Also assume that none of the $q_i$ can be written as
$t^m f(t^7)$ for some $f(t) \in \qq(\zeta)[t,t^{-1}]$ and integer $m$.

Suppose  that $k_1,k_2,\ldots,k_n$ are non-negative integers and that 
$a_{ij}\in\zz_7$, $i=1,2,\ldots,n,$ \ $j=1,2,\ldots, k_i$, and let 
$$R(t)=\left(\prod_{j=1}^{k_1}q_1(\zeta^{a_{1j}}t)\right)
\cdots \left(\prod_{j=1}^{k_n}q_n(\zeta^{a_{nj}}t)\right).$$ Then
$R(t)$ is a norm in $\qq(\zeta)[t,t^{-1}]$ if and only if
 each $a\in\zz_7$ appears an even number of times in the vector
$(a_{\ell 1},a_{\ell 2},\ldots,a_{\ell k_\ell })$ for each
$\ell\in\{1,2,\ldots,n\}$.

\end{lemma}
\begin{proof} Each irreducible factor of $R(t)$ is of the form
$q_i(\zeta^j t)$ for some $i$ and $j$, since the irreducibility of a
polynomial
$q(t)$ implies the irreducibility of $q(\zeta^a t)$.  Each of these
factors, $q_i(\zeta^j t)$, satisfies
$\overline{q_i(\zeta^j t)} =q_i(\zeta^j t)$.  Hence, if $R(t) =
Q(t)\overline{Q(t)}$ for some
$Q(t)$, then the exponent of each irreducible factor of $R(t)$ must be
even.

It remains to observe that the set of polynomials $\left\{ q_i(\zeta^j
t)
\right\}_{i=1,\ldots ,n ,\ j= 0, \ldots ,6}$ are distinct.  Suppose that
$q_\alpha(\zeta^r t) = q_\beta(\zeta^s t)$.  Then a change of variable
shows that $q_\alpha(t) = q_\beta(\zeta^{s-r} t)$. Since the $q$ are
integer polynomials, the only way these can be equal is if $s-r = 0$. 
(This is where we use the fact that the $q_i(t)$ are not polynomials in
$t^7$.)  But this results in the equation $q_\alpha(t) = q_\beta(t)$,
implying that $\alpha = \beta$.
\end{proof}

\subsection{Case I: Metabolizers with odd characters}  We turn now to
the proof of Theorem \ref{phi1ker}, that is, to showing that
the kernel of
$\phi_1$ contains $\zz^\infty$. There are basically two separate
cases in the argument, depending on whether or not the space of
characters that vanish on a given metabolizer contains an {\sl odd
vector}, to be defined below.

Let ${\cal S}=\{S_1,S_2,\cdots\} $ be an infinite collection of isotopy
classes of oriented knots in
$S^3$ so that their Alexander polynomials $\Delta_{S_i}(t)$ satisfy
\begin{enumerate}
\item $\Delta_{S_i}(t)$ is irreducible in
$\qq(\zeta)[t,t^{-1}]$.

\item$\Delta_{S_i}(t)$ is relatively prime to
$\Delta_{S_j}(t)$  for all
$i,j$.

\item  $\Delta_{S_i}(t)$ is not of the form $t^m f(t^7)$ for some $m$
and rational polynomial $f$.
\item $\Delta_{S_i}(\zeta^a t)$ is relatively prime to the 49 elements
$\delta(K;\chi),\chi\in H^1(M_3;\zz_7)$.

\end{enumerate}

Such a collection is easily constructed.  Consider for instance the
polynomial  $S(t) = at^2 - (2a+1)t - a$, $a > 0$.  If $a$ is not of the
form
$r(r+1)$ for some integer $r$, then $S(t)$ is irreducible and has real
roots
${(2a+1) \pm \sqrt{4a+1} \over 2a}$.  If $S(t)$ were reducible  over
$\qq(\zeta_7)$, then we would have the extension $\qq
\subset \qq(\sqrt{4a+1}) \subset \qq(\zeta_7)$.  However, $\qq(\zeta_7)$
is a degree 6 Galois extension of $\qq$ and hence  contains a unique
degree  2 extension of $\qq$.  As we show next, that extension is not
real while $\qq(\sqrt{4a+1}) $ is, leading to the desired contradiction.

To see that the degree 2 extension $F$ of $\qq$ contained in 
$\qq(\zeta_7)$ is complex, proceed as follows. We have that $F$ is the
fixed set of the order 3 Galois automorphism $\sigma$ of
$\qq(\zeta_7)$,  defined by
$\sigma(\zeta_7) = \zeta_7^2$.  A fixed element of $\sigma$ is $\alpha =
\zeta_7 + \zeta_7^2 +\zeta_7^4$.  One checks that $\alpha$ satisfies $t^2
+t + 2 = 0$, so $\alpha = {-1 + \sqrt{-7} \over 2}$ and $F =
\qq(\sqrt{-7})$. 

So, to construct the $S_i$ pick an infinite sequence of integers $a$ as
above.  The resulting polynomials will satisfy the first three
conditions and all but a finite number of them satisfy the last
condition.  These polynomials occur as the Alexander polynomials of
twisted doubles of the unknot.

Next suppose we are given
\begin{enumerate}
\item  a positive number $n$,
\item   a choice of knots $J_i$, $i=1,\ldots ,n$ such that each
$J_i\in \cals$, and
\item a choice of signs $\epsilon_i\in\{\pm 1\}$ for
$i=1,\ldots , n$, such that if $J_i=J_j$ then $\epsilon_i=\ep_j$.
\end{enumerate}

Denote by
$L$ the connected sum of the $n$ (oriented) knots $$L=\ep_1K^*_{J_1}
\# \ep_2K^*_{J_2}\#
\cdots \#\ep_nK^*_{J_n}$$ and let $M_3(L)$ denote 3--fold branched
cover of $L$. Notice that performing $n$ mutations transforms $L$
into the slice knot $\ep_1K_{J_1}
\# \ep_2K_{J_2}\#
\cdots \#\ep_nK_{J_n}$, and hence $L$ lies in the kernel of
$\phi_1\!:\calc\to\calm$.

We will show that
$L$ is not slice by analyzing characters $\chi$ in
$H^1(M_3(L);\zz_7)$. This will show that the kernel of $\phi_1$
contains an infinite direct sum of $\zz$.

The 3--fold branched cover of $S^3$ branched over $L$, $M_3(L)$, is the
connected sum of
$n$ copies of $M_3(\ep_iK_{J_i}^*)$, since the branched cover of $S^2$
branched over 2 points is a 2--sphere. Thus
$$H^1(M_3(L);\zz_7)=\oplus_{i=1}^n
H^1(M_3(\ep_iK_{J_i}^*);\zz_7).\eqno{(1)}$$
At this point we have identified $H_1(M_3(K))$ with
$H_1(M_3(K_J^*))$. We also need to identify $H_1(M_3(K_J^*))$ with
$H_1(M_3(-K_J^*))$, along with their $\zz_3$ structures. First, since
$-(S^3,K)=(-S^3,-K)$, $M_3(K_J^*) = -M_3(-K_J^*)$. Thus their first
homologies are naturally identified. Next, we need to check that the
$\zz_3$ actions agree under this identification. This follows from the
fact that an oriented curve in
$S^3 - K$ that links $K$ algebraically once also links
$-K$ algebraically once in $-S^3$. In particular, $\delta(K_J^*,\chi) =
\delta(-K_J^*,\chi)$ for any $J$.

The covering transformation acts diagonally in the sum of Equation (1),
hence we can write $$H^1(M_3(L);\zz_7)=\oplus_{i=1}^n (C_2\oplus
C_4),$$ and characters $\chi$ in
$H^1(M_3(L);\zz_7)$ will be expressed as an $n$--tuple
$$\chi=(a_1\chi_2+b_1\chi_4, a_2\chi_2+b_2\chi_4,\cdots,
a_n\chi_2+b_n\chi_4).\eqno{(2)}$$
\begin{lemma} \label{discrform}For the character
$\chi\in H^1(M_3(L);\zz_7)$ of Equation (2), the\break Casson--Gordon
discriminant
$\delta(L;\chi)$ equals

\begin{eqnarray*}
\prod_{i=1}^n  \mbox{\rm \LARGE
 \raisebox{-.02in}{(}}  \delta_i \cdot\Delta_{J_i}(\zeta^{a_i+b_i}t)
&&\hskip-.26in\Delta_{J_i}(\zeta^{2a_i+4b_i}t)
\Delta_{J_i}(\zeta^{4a_i+2b_i}t)   \\
 &&\Delta_{J_i}(\zeta^{-a_i-b_i}t)
\Delta_{J_i}(\zeta^{-2a_i-4b_i}t)\Delta_{J_i}(\zeta^{-4a_i-2b_i}t)
\mbox{\rm \LARGE
 \raisebox{-.02in}{)}}  
\end{eqnarray*} where
$\Delta_{J_i}(t)$ denotes the Alexander polynomial of the knot
${J_i}$ and $\delta_i=\delta(K,{a_i\chi_2+b_i\chi_4})$.
\end{lemma}
\begin{proof} Using Theorem \ref{cohomology}  one sees that since the
character
$a_i\chi_2+b_i\chi_4$ takes the value
$a_i+b_i$ on the lift $\tilde{B}_1$ of $B_1$, it takes the value
$2a_i+4b_i$ on the translate
$T\tilde{B}_1$ and takes the value $ 4a_i+2b_i$ on $T^2\tilde{B}_1$.

Similarly
$a_i\chi_2+b_i\chi_4$ takes the value $-a_i-b_i$ on $\tilde{B}_2^*$
 and hence takes the value $-2a_i-4b_i$ on
$T\tilde{B}_2^*$ and $-4a_i-2b_i$ on
$T^2\tilde{B}_2^*$.

The lemma now follows from Theorem \ref{deltathm}.\end{proof} 

To show that $L$ is not slice, we will prove that for each invariant
metabolizer
$A\subset H_1(M_3(L);\zz)$ there exists a $\zz_7$ character $\chi$
vanishing on $A$ for which $\delta(L,\chi)\in \qq(\zeta)[t,t^{-1}]$ is
not a norm.

Given a metabolizer $A
\subset H_1(M_3(L))$ which is invariant under the $\zz_3$ action
induced by the covering transformation, let $A^* \in H^1(M_3(L);\zz_7)$
denote the subgroup of $\zz_7$ characters vanishing on $A$ (we identify
$H^1(M_3(L);\zz_7)$ with Hom$(H_1(M_3(L));\zz_7)$). Since
$A$ is $\zz_3$ invariant, so is $A^*$, hence $A^*$ is spanned by
eigenvectors of the
$\zz_3$ action. Hence $A^*$ decomposes into $2$ and $4$ eigenspaces (as
before, the $1$ eigenspace is trivial since $H^1(S^3;\zz_7)=0$):
$$A^*=A^*_2\oplus A^*_4.$$
For $i=2$ or $4$, any character $\chi$ in $A_i^*$ can be written in the
form
$\chi=(a_1\chi_i,a_2\chi_i,\cdots,a_n\chi_i)$. 

\begin{defn} An eigen-character
$\chi=(a_1\chi_i,a_2\chi_i,\cdots,a_n\chi_i)$ in $A^*_i$ is called
 {\em odd }  if an odd number of the coefficients $a_j$ are
non-zero.\end{defn}

\begin{lemma}\label{odd} If $\chi\in A^*_i$ is odd, then
$\delta(L,\chi)$ is not a norm in
$\qq(\zeta)[t,t^{-1}]$.\end{lemma}
\begin{proof} Given $a\in \zz_7$, the set
$$\{\zeta^a,\zeta^{2a},\zeta^{4a},
\zeta^{-a},\zeta^{-2a},\zeta^{-4a}\}$$ equals
 $\{\zeta,\zeta^{2},\zeta^{3},
\zeta^{4},\zeta^{5},\zeta^{6}\}$  if $a\ne 0$, and equals
 $\{1,1,1,1,1,1\}$  if $a=0$. 

Suppose
$\chi=(a_1\chi_2,a_2\chi_2,\cdots,a_n\chi_2)\in A^*_2$. Then by Lemma
\ref{discrform} and the choice of $J$
$$\delta(L,\chi) =\delta\!\prod_{i|a_i=0}\!\Delta_{J_i}(t)^6
\prod_{i|a_i\ne 0}\Delta_{J_i}(\zeta t) \Delta_{J_i}(\zeta^2
t)\Delta_{J_i}(\zeta^3 t)\Delta_{J_i}(\zeta^4 t)\Delta_{J_i}(\zeta^5
t)\Delta_{J_i}(\zeta^6 t)$$ Where $\delta=\prod_i\delta(K,a_i\chi)$.

By Lemma \ref{squares} $\Delta_{J_i}(t)^6$ is  a norm. Using the
hypotheses on
$\Delta_{J_i}$, Lemma \ref{squares}, and the fact that $\chi$ is odd,
we see that, modulo norms, $\delta(L,\chi)$ equals a product  
$$\delta \cdot \prod_{r\in B} \Delta_{S_r}(\zeta t)
\Delta_{S_r}(\zeta^2 t)\Delta_{S_r}(\zeta^3 t)\Delta_{S_r}(\zeta^4
t)\Delta_{S_r}(\zeta^5 t)\Delta_{S_r}(\zeta^6 t) \eqno{(1)}$$
 where
$B\subset\{1,2,\ldots\}$ is the set of those $r$ so that the set
$\{i|J_i=S_r \hbox{ and } a_i\ne0\}$ contains an odd number of elements.

 Since $\chi$ is odd,
$B$ is non-empty  and so by Lemma \ref{squares} the expression in
Equation (1) is not a norm.

The same argument applies to $A_4^*$.\end{proof}

\begin{cor}\label{oddcase} Given a metabolizer $A$, if one of $A_2$ or
$A_4$ contains an odd character $\chi$ then $\delta(L,\chi)$ is not a
norm .\qed\end{cor}

In preparation for the next lemma, notice that $A^*\subset
H^1(M_3(L);\zz_7)=(\zz_7)^{2n}$ always has dimension at least $n$. This
follows from the fact that the metabolizer $A\subset 
H_1(M_3(L))=(\zz_{49})^{2n}$ has order
$\sqrt{|(\zz_{49})^{2n}|}= 49^n$, and so
$$A^*=\hbox{Hom}(H_1(M_3(L))/A,
\zz_7)=\hbox{Hom}((\zz_{49})^{2n} /A;\zz_7)$$ has dimension at least
$n$.

\begin{lemma}\label{existodd} If one of $A_2^*$ or $A^*_4$ has
dimension greater than
$n/2$, then it contains an odd character. \end{lemma} 

\begin{proof} Suppose that $A_\ell^*$ has dimension $k>  n/2$. Choose some
basis for
$A_{\ell}^*$. Each basis element is written
$a^i=(a^i_1\chi_{\ell}, a^i_2\chi_{\ell}, \cdots ,a^i_n\chi_{\ell})$
which we consider as the $i$th row of a $k\times n$ matrix $M$. 

Elementary row operations, and interchanging of columns, which does not
affect whether a row is odd (ie, has an odd number of non-zero entries)
transforms the matrix $M$ to a matrix of the form
$\left(I\ B\right)$, where $I$ is a $k\times k $ identity matrix and
$B$ is a $k\times(n-k)$ matrix. If some row of $B$ is not odd, then the
corresponding row of $(I\ B)$ is odd and we have the desired odd
character.

Thus assume every row of $B$ is odd. Since
$B$ has more rows than columns, some non-trivial linear combination of
the rows of $B$ gives the zero vector, say $\sum_i s_i\
\hbox{row}_i(B)=0$. If an odd number of the $s_i$ are non-zero, then
$\sum_i s_i\ \hbox{row}_i(I\ B) $ is odd.

Finally, if an even number of the $s_i$ are nonzero, choose $j$ so that
$s_j\ne 0$, and pick $f\in\zz_7 $ non-zero and different from
$s_j$. Then $-f\cdot\hbox{row}_j(I\ B) +\sum_i s_i\ 
\hbox{row}_i(I\ B)
$ is odd since $f\cdot $row$_j(B)$ is odd. \end{proof}

\noindent{\bf Remark}\qua Notice that the proof also works if
$A^*_{\ell}$ has dimension
$n/2$ and the square matrix $B$ is singular. We will use this extension
below. 

\medskip

If $n$ is odd, then for each $A$ one of $A^*_2$ and
$A^*_4$ has dimension greater than $n/2$, and hence contains an odd
character. This shows
$nK_J^*$ is not slice for $n$ odd or, more generally, that $L$ is not
slice if $n$ is odd. But when
$n$ is even there are some remaining cases when neither $A_i^*$
contains an odd vector. For example, consider the situation when
$L=K^*_{S_1}\#K_{S_1}^*$ and
$$A^*_2=\hbox{span}\{ (\chi_2, a\chi_2)
\}$$ and
$$A^*_4=\hbox{span}\{ (\chi_4, b\chi_4)\}$$ with $a,b$ non-zero. Then
$A_2^*$ and $A_4^*$ do not contain any odd vectors. Using Lemma
\ref{discrform} one computes that with
$\chi=(\chi_2, a\chi_2)+(\chi_4, - a^{-1}\chi_4)$,
$$\delta(2K_J^*,\chi)=\delta(K,\chi_2+\chi_4)\delta(K,a\chi_2-a^{-1}\chi_4)
\cdot
\Delta_J(\zeta^2t)\Delta_J(\zeta^5 t)$$ modulo norms, which is
non-trivial by Lemma \ref{squares}. The choice $b=-a^{-1}$ is crucial,
and to establish this equality (appropriately generalized) we will have
to switch to
$\zz_{49}$ coefficients and examine the linking form more carefully.

\subsection{Case II: Metabolizers with no odd characters
}\label{kerphi1}

The proof that
$L$ is not slice  is reduced to treating the case when $A^*_2$ and
$A^*_4$ both have dimension $n/2$ (in particular $n$ is even) and
neither $A_2^*$ nor
$A_4^*$ contains an odd character. We will ultimately show that in this
case there nevertheless exists a $\chi\in A^*$ such that
$\delta(L,\chi)$ is not a norm. The analysis turns out to be  more
complicated and, in particular, will require the use of the linking
form on $M_3$ and $\zz_{49}$ coefficients.

The linking form on the 3--fold branched cover of a knot $Z$ is a
non-singular symmetric pairing
$$ lk\!:H_1(M_3(Z);\zz)\times H_1(M_3(Z);\zz)\to \qq/\zz$$ Taking $Z=K$
or
$Z=L$ or $Z=K_J^*$, we see that the linking form takes its values in
$\zz_{49}\subset \qq/\zz$. Hence the adjoint of the linking form
defines an isomorphism $$\ell\!:H_1(M_3(Z);\zz)\to H^1(M_3(Z);\zz_{49})
$$ by $\ell(x)(y)=lk(x,y)$. This in turn defines a pairing
$$lk^*\!:H^1(M_3(Z);\zz_{49})\times H^1(M_3(Z);\zz_{49})\to \zz_{49}$$ by
$lk^*(x,y)=lk(\ell^{-1}(x),\ell^{-1}(y))$. 

\begin{lemma}\label{lkform} In the basis
$\chi_{18},\chi_{30}$ for $H^1(M_3(K);\zz_{49})$, the form
$ lk^*$ has matrix
$$\begin{pmatrix} 0 & u \\ u & 0 \end{pmatrix}$$ for some unit
$u\in\zz_{49}$.
\end{lemma}
\begin{proof} Since $T(\ell(x))=\ell(T^{-1}(x))$ for $x\in
H_1(M_3(K);\zz)$, it follows that
$T$ takes $18$ eigenvectors to $30$ eigenvectors and vice-versa. Hence
$ lk^*(\chi_{18},\chi_{18})= 18^{-2}lk^*(\chi_{18},\chi_{18})$ and so
$lk^*(\chi_{18},\chi_{18})=0$ since $18^{-2}-1=17=26^{-1}$. Similarly
$lk^*(\chi_{18},\chi_{18})=0$. The off--diagonal entries are equal since
the form is symmetric, and $u$ must be a unit since the form is
non-singular.\end{proof}

To an invariant metabolizer $A\subset
H_1(M_3(L);\zz)\cong(\zz_{49})^{2n}$ we have associated the group
$A^*\subset H^1(M_3(L);\zz_7)\cong(\zz_{7})^{2n}$ of $\zz_7$ characters
which vanish on $A$.

\begin{lemma} \label{linking}Let $A\subset (\zz_{49})^{2n}$ be an
invariant metabolizer. Suppose that
$A^*$ is isomorphic to $(\zz_7)^{n}$. Suppose that $a_1\cdots,a_n$ and
$ b_1,\cdots, b_n$ are integers so that $(a_1\chi_2,\cdots,
a_n\chi_2)\in A^*_2$ and $(b_1\chi_4,\cdots ,b_n\chi_4)\in A^*_4$. Then
$\sum_i \ep_i a_i b_i=0\pmod 7$.\end{lemma}

\proof Write $M$ for $M_3(L)$. We always identify
the group $H^1(M;\zz_{7})$ with
$7H^1(M;\zz_{49})\subset H^1(M;\zz_{49})$. Recall that
$\chi_2=7\chi_{30}$ and $\chi_4= 7\chi_{18}$ and therefore we have 
$$(a_1\chi_2,\cdots, a_n\chi_2) = 7 (a_1\chi_{30},\cdots, a_n\chi_{30})$$
and $$(b_1\chi_4,\cdots ,b_n\chi_4)=7 (b_1\chi_{18},\cdots
,b_n\chi_{18}).$$
Since $A$ is a subgroup of
$(\zz_{49})^{2n}$, $A\cong\zz_{49}^a\oplus \zz_7^b$. Since $A$ is a
metabolizer, its order is $49^{n}$ and so $2a+b=2n$. The group $A^*$
equals Hom$(H_1(M)/A;\zz_7)$. But
$H_1(M)/A\cong A$. To see this, note that the sequence
$$0\to A\to H_1(M)\mapright{\ell}
\hbox{Hom}(A;\qq/\zz)\to 0$$ is exact, where $\ell(x)(y)=lk(x,y)$,
since the linking form is non-degenerate. But since $A$ is a finite
abelian group,
$\hbox{Hom}(A;\qq/\zz)$ is isomorphic to $A$. 

Thus $(\zz_7)^{n}\cong A^*\cong$ Hom$(A,\zz_7)\cong (\zz_7)^{a+b}$, so
$n=a+b$. Together with $2a+b=2n$ this implies that $a=n$ and
$b=0$, so $A\cong (\zz_{49})^{n}$. 

Now let $\tilde{A}\subset H^1(M;\zz_{49})$ denote the set of
$\zz_{49}$ characters which vanish on $A$; the map $\ell$ takes $A$
isomorphically to $\tilde{A}$, so that
$\tilde{A}\cong (\zz_{49})^n$. It follows easily that $A^*=7
\tilde{A}$. 

Therefore, there exist $x,y\in H^1(M;\zz_7)$ so that $$
(a_1\chi_{30},\cdots, a_n\chi_{30})+7x $$ and
$$(b_1\chi_{18},\cdots ,b_n\chi_{18}) + 7y $$ are both in $\tilde A$.
Since $M$ is a
connected sum, the linking form
$lk^*$ splits according to the splitting of $H^1(M;\zz_{49})=
\sum_{i=1}^n H^1(M_3(\ep_iK_{J_i}^*);\zz_{29})$ . On each summand, the linking matrix is $$\ep_i
\begin{pmatrix} 0 & u \\ u & 0 \end{pmatrix},$$  using Lemma
\ref{lkform} and the fact that
$M_3(-Z)=-M_3(Z)$.

Thus
$$0 =lk\left((a_1\chi_{30},\cdots,
a_n\chi_{30})+7x,(b_1\chi_{18},\cdots ,b_n\chi_{18}) + 7y\right)$$
$$=u\sum_i\ep_i a_ib_i\pmod 7.\eqno{\qed}$$

We can now return to the argument that $L$ is not slice. We have treated
all cases except when $n=2k$,
$A_2^*$ and
$A_4^*$ are both
$k$ dimensional and contain no odd eigencharacters. 

A vector in $A_2^*$ can be expressed as an $n$--tuple in $\zz_7$ using
the correspondence
$$(a_1,\cdots,a_n)\leftrightarrow(a_1\chi_2,\cdots a_n\chi_2).$$ This way
any basis for
$A_2^*$ determines a $k\times 2k $ matrix $(D\ E )$ whose rows are the
basis vectors. Similarly a basis for $A_4^*$ determines a $k\times 2k$
matrix $(F\ G)$. 

As remarked after Lemma \ref{existodd}, if the matrix $(D\ E)$ can be
reduced to
$(I\ B)$ with $B$ singular by performing elementary row operations and
column interchanges, then $A^*_2$ contains an odd character
$\chi$ and so by Corollary \ref{oddcase}     for this odd character
$\delta(L;\chi)$ is not a norm. Thus it remains only to treat the case
when  in any such reduction the matrix $B$ is non-singular, and
similarly for the 4-eigenspace  
$A^*_4$.

Lemma \ref{linking} above implies that
$$\begin{pmatrix} D & E \\ \end{pmatrix} S \begin{pmatrix} F^t \\ G^t
\end{pmatrix}=0$$ where $S$ is the diagonal matrix with diagonal
entries the signs
$\ep_1,\cdots, \ep_n$. 

There exists an invertible matrix $K\in GL_k(\zz)$ and a
$2k\times 2k$ permutation matrix $P$ so that $K(D\ E) P^{-1}=(I\ B)$,
since elementary row operations and column interchanges reduce the rank
$k$ matrix $(D\ E)$ to $(I\ B)$ with $B$ non-singular. Changing the
basis of $A_2^*$ we may assume that $K=I$, so that $(D\ E)=(I\ B)P$.  
 
Define $k\times k$ matrices   $X$ and $Y$ by $(X\ Y)=  (F\ G)P^{-1}$,
and define two $k\times k$ diagonal matrices $S_1$ and $S_2$ with
diagonal entries the signs $\epsilon_i$ permuted by $P$ by
$$\begin{pmatrix}S_1 & 0\\ 0 &S_2\end{pmatrix}=PSP^{-1}.$$
Then
\begin{eqnarray*} 0&=&\begin{pmatrix} D&E \end{pmatrix} S
\begin{pmatrix} F^t\\ G^t \end{pmatrix}\\ &=&\begin{pmatrix} I&B
\end{pmatrix}\begin{pmatrix}   S_1& 0\\ 0&S_2 \end{pmatrix}
\begin{pmatrix} X^t\\ Y^t \end{pmatrix}\\
&=&S_1X^t+BS_2Y^t\end{eqnarray*}

Hence $Y=-XS_1(B^t)^{-1}S_2^{-1}$ and so
$$(F\ G)=(X\ Y) P=X(I\ -S_1(B^t)^{-1}S_2)P.$$ The matrix $X$ must
therefore be invertible and so by changing basis of
$A_4^*$ by the matrix $XS_1$ we can assume
$(F\ G)= (I\ -(B^t)^{-1})PS.$

By assumption $A^*_2$ and $A^*_4$ contain no odd vectors.  Thus we are
poised to apply the following lemma, the proof of which is given later in
this section.

\begin{lemma}\label{larsen} Let $E$ be a non-singular $k \times k$
matrix over
$\zz_{p}$ for a prime
$p>2$. Suppose that the subspace of $(\zz_p)^{2k}$ spanned by the rows
of the $k\times 2k$ matrix $(I\ E)$ contains no odd vectors (ie, every
vector in this span has an even number of non-zero entries). Then $E$
is obtained from a diagonal matrix by permuting the columns.
\end{lemma}

The proof that $L$ is not slice is completed as follows. Lemma
\ref{larsen} implies that the matrix $B$ is obtained from a diagonal
matrix by permuting the columns, say $B=CQ$ for a diagonal matrix
$C=$diag$(c_1,\ldots,c_n)$ and a $k\times k$ permutation matrix $Q$.
Write 
$R=\begin{pmatrix}I& 0\\ 0 &Q\end{pmatrix}P$; thus $R$ is a $2k\times
2k$ permutation matrix and
$$(D\ E)= (I\ C)R,\ \ (F\ G)= (I\ -C^{-1}  )RS.$$
By reordering the summands of $L=\ep_1 K_{J_1}^*\#\cdots\#
\ep_n K_{J_n}^*$ we may assume that $R$ is the identity permutation. 
Thus $A^*_2$ contains the character  
$$\alpha_2=( \chi_2,0,\ldots,0,c_1\chi_2,0,\ldots,0)$$ and $A^*_4$
contains the character $$\alpha_4=(
\ep_1\chi_4,0,\ldots,0,-\ep_{k+1}c_1^{-1}\chi_4,0,\ldots,0)$$ where for
each vector the second non-zero entry   is in the $k+1$st position. 
Thus any linear combination of these characters lies in 
$A^*=A_2^*\oplus A_4^*$.

Let
$\chi=\alpha_2 +\ep_1\alpha_4$. We claim that $\delta(L,\chi)$ is not a
norm.

For ease of notation write $J=J_1$ and $J'=J_{k+1}$,
$\ep=\ep_1\ep_{k+1}$,  
$c=c_1$, and $\bar{c}=c^{-1}\pmod{7}$. 

Lemmas
\ref{discrform}  and \ref{squares} show that, modulo norms,
$\delta(L,\chi)$ is the product of the four terms:
\begin{enumerate}  
\item $\delta(K,\chi_2+\ep\chi_4)$, 
\item  $\delta(K,c\chi_2 -\ep'c^{-1}\chi_4)$,
 \item $\Delta_{J}(\zeta^{2}t)\Delta_{J}(\zeta^{6}t)
\Delta_{J}(\zeta^{6}t)\Delta_{J}(\zeta^{5}t)
\Delta_{J}(\zeta t)
\Delta_{J}(\zeta t)\equiv\Delta_{J}(\zeta^{2}t)\Delta_{J}(\zeta^{5}t)$
(modulo norms)
\item
$\Delta_{J'}(\zeta^{c-\ep\bar{c}}t)
\Delta_{J'}(\zeta^{2c-4\ep\bar{c}}t)
\Delta_{J'}(\zeta^{4c-2\ep \bar{c}}t) \cdot$ 

$\hskip1.5in \Delta_{J'}(\zeta^{-c+\ep \bar{c}}t)
\Delta_{J'}(\zeta^{-2c+4\ep\bar{c}}t)
\Delta_{J'}(\zeta^{-4c+2\ep\bar{c}}t)$.\end{enumerate}

There are two cases to consider, depending on whether or not $J=J'$.
 If $J=J'$, ie, $J_{1}\ne J_{k+1}$, then 
$\delta(L,\chi)$ is not a square since the third term above is not a
norm and cannot cancel with any of the other three terms by the choice 
of the $S_r$.

If $J_1=J_{k+1}$, then $\ep_1=\ep_{k+1}$, so that $\ep=1$. Then the
exponents of $\zeta$ in the fourth term above are
$$\{ c-\bar{c},
2c-4\bar{c},4c-2\bar{c},-c+\bar{c},-2c+4\bar{c},-4c+2\bar{c}\}.$$ The
reader can check that for each non-zero choice of $c$ in $\zz_7$, this
(unordered) set is just $\{0,0,2,2,5,5\}$, so that the fourth term is
always a norm. Since the third term is not a norm, the product of the
third and fourth terms is not a norm, and so $\delta(L,\chi)$ is not a
norm.

This completes the argument that $L$ is not slice, since for each
invariant metabolizer $A$ we have shown how to construct a character
$\chi\in A^*$ for which $\delta(K,\chi)$ is not a norm.  Hence Theorem
\ref{phi1ker} is proven.

\subsection{The proof of Lemma \ref{larsen}}

\begin{proof} Consider the columns of $E$ as a basis $e_i$ of
$V=(\zz_p)^k$. Let
$b_1,\cdots, b_n$ denote the standard basis of $V$. Let
$V^*=$Hom$(V,\zz_p)$. The condition that every linear combination of
the rows of
$(I\ E)$ is even can be restated by saying that for each $x\in V^*$,
the sets
$$\{ j| x(e_j)\ne 0\}\ \hbox{and}\ \{j|x(b_j)\ne 0\}$$ have the same
cardinality modulo
$2$.

Another way to state this is as follows. Let $h_j\!:V^*\to \cc$ be the
function 
$$h_j(x) =\begin{cases} \ \ 1& \hbox{if}\ x(e_j)=0,\\ -1& \hbox{if}\
x(e_j)\ne 0
\end{cases}$$ and let $H\!:V^*\to \cc$ be the product
$H(x)=\prod_{j=1}^k h_j(x)$. Similarly define $g_j\!:V^*\to \cc$ by
$$g_j(x)=\begin{cases}
\ \ 1& \hbox{if}\ x(b_j)=0,\\ -1&\hbox{if}\ x(b_j)\ne 0 \end{cases}$$
and $G\!:V^*\to \cc$ be the product $G(x)=\prod_{j=1}^k g_j(x)$. Then the
hypothesis is equivalent to the statement that $H=G$. 

Let $\hat{H}\!:V\to \cc$ be the Fourier transform $$\hat{H}(v)=
\sum_{x\in V^*}H(x)e^{{{2\pi i}\over p}x(v)}.$$ Define $\hat{G}\!:V\to
\cc$ similarly using $G$.

We compute $\hat{H}$ as follows. Let $\{e_j^*\}\subset V^*$ denote the
dual basis to the
$e_j$. Then $$\hat{H}(v)= \sum_{x_1\in\zz_p}\cdots\sum_{x_k\in\zz_p}
H(\Sigma_i x_i e_i^*)e^{{{2\pi i}\over
p}\Sigma_ix_ie_i^*(v)}.\eqno{(3)}$$
From its definition one sees that $h_j(\sum_i x_i e_i^*)=h_j(x_j 
e_j^*)$, so that
$$H(\Sigma_i x_i e_i^*)=\prod_j h_j(x_je_j^*).$$ Substituting this into
Equation (3) and simplifying one gets $$\hat{H}(v)
=\prod_j\left(\sum_{x\in \zz_p}h_j(xe_j^*) e^{{{2\pi i}\over p} x
e_j^*(v)}\right).$$
Now
$$ \sum_{x\in \zz_p}h_j(xe_j^*) e^{{{2\pi i}\over p} x e_j^*(v)} =
1-\sum_{x\in \zz_p-0} e^{{{2\pi i}\over p} x e_j^*(v)}.$$
Notice that 
\[
\sum_{x\in \zz_p-0} e^{{{2\pi i}\over p} x e_j^*(v)} 
=\begin{cases}
p-1&\text{  if } e^*_j(v)=0,\\
-1&\text{  if } e^*_j(v)\ne 0.
\end{cases}\]
 Thus for each $v$,
$\hat{H}(v)= 2^{\alpha(v)}(2-p)^{k-\alpha(v)}$ where $\alpha(v)$ is the
cardinality of the set $\{j| e_j^*(v)\ne 0\}$. Similarly
$\hat{G}(v)= 2^{\beta(v)}(2-p)^{k-\beta(v)}$ where $\beta(v)$ is the
cardinality of the set $\{j| b_j^*(v)\ne 0\}$. Since $H=G$, $\hat{H}
=\hat{G} $, and since $p> 2$,   
$$2^{\alpha(v)-\beta(v)}=(2-p)^{\alpha(v)-\beta(v)}.$$
This implies that
 $\alpha(v)=\beta(v)$  since $p\ne 4$  (recall that $p$ is prime).  Now
$\alpha(e_j)=1$ and therefore $\beta(e_j)=1$. Hence
 for each $j$ there exists a
unique $k_j$ so that $ b_{k_j}^*(e_j)\ne 0$. Since $E$ is non-singular
it follows that $j\mapsto k_j$ is a permutation, and that
$e_j=(b_{k_j}^*(e_j))b_{k_j}$. In other words, $E$ is obtained from a
diagonal matrix by permuting the columns.\end{proof}

\end{document}